\documentclass[smallextended,referee,envcountsect]{article}
\usepackage[fontset=ubuntu, nofonts]{ctex}
\usepackage{amsfonts,amsmath,amssymb}
\usepackage{newtxtext, newtxmath}
\usepackage{ctex}
\usepackage{float}
\usepackage{latexsym}
\usepackage{mathrsfs}

\usepackage[pdftex]{graphicx}
\usepackage[ruled]{algorithm2e}
\usepackage{algpseudocode}
\usepackage{indentfirst}
\usepackage{tabularx}
\usepackage{booktabs}
\usepackage{longtable}
\usepackage{multirow}
\usepackage{caption}
\usepackage{subeqnarray}
\usepackage{verbatim}
\usepackage{tikz}
\usepackage[english]{babel}
\usepackage{siunitx}
\usepackage{subcaption}
\usepackage{booktabs}
\usepackage{graphicx}
\usepackage{array}
\usepackage{empheq}
\usepackage{xcolor}

\newcommand{\rothead}[1]{\rotatebox{45}{\parbox{1.8cm}{\centering #1}}}
\newtheorem{theorem}{Theorem}[section]
\newtheorem{lemma}[theorem]{Lemma}
\newtheorem{assumption}[theorem]{Assumption}

\newtheorem{definition}[theorem]{Definition}

\newtheorem{remark}[theorem]{Remark}
\numberwithin{equation}{section}
\newenvironment{proof}[1][Proof]{\textbf{\varphi1.} }
{\ \rule{0.75em}{0.75em}\smallskip}

\usepackage[pdftex,unicode,colorlinks,
linkcolor=blue]{hyperref}

\def\lista\varphi1
{{ \itemindent 0.0cm \labelsep .2cm \leftmargin 0.8cm \rightmargin
		0.0cm \labelwidth 0.6cm \topsep 0.0mm
		\parsep 0.0mm
		\itemsep 0.0mm
		\begin{list}{}
			{ \setlength{\leftmargin}{.8cm} \setlength{\rightmargin}{0.0cm}
				\setlength{\parsep}{0.0mm} \setlength{\topsep}{.0mm}
				\setlength{\parskip}{.0cm} \setlength{\itemsep}{.0cm} }
			{\varphi1}\end{list}} }

\begin{document}
	
	\title{\Large \bf
		A Proximal Stochastic  Gradient Method with Adaptive Step Size and Variance Reduction for Convex Composite Optimization}

	\author{
		Changjie Fang \,\footnote{\, School of Mathematics and Statistics, Chongqing University of Posts and Telecommunications,	Chongqing 400065, China.
			E-mail: {\tt fangcj@cqupt.edu.cn.}},
		\ \
		Hao Yang \,\footnote{\, School of Mathematics and Statistics, Chongqing University of Posts and Telecommunications, Chongqing 400065, China.
			E-mail: {\tt yanghao55255@163.com.}},
		\ \
		Shenglan Chen \,\footnote{\, School of Mathematics and Statistics, Chongqing University of Posts and Telecommunications,	Chongqing 400065, China.
			E-mail: {\tt chensl@cqupt.edu.cn.}}
	}

	\date{}
	\maketitle

	\noindent {\bf Abstract}
	In this paper, {we propose a proximal stochastic gradient} algorithm (PSGA) for solving composite optimization problems by incorporating variance reduction techniques and an adaptive step-size strategy. In the PSGA method, the objective function consists of two components: one is a smooth convex function, and the other is a non-smooth convex function. We establish the strong convergence of the proposed method, provided that the smooth convex function is $L$-smooth. We also prove that the expected value of the error between the estimated gradient and the actual gradient converges to zero. Furthermore, we get an \( O(\sqrt{1/k}) \) convergence rate for our method. Finally, the effectiveness of the proposed method is validated through numerical experiments on Logistic regression and Lasso regression.
	\vspace{5mm}	
	
	\vspace{5mm}

\section{Introduction}
\noindent
{In this paper, we focus on large-scale convex composite optimization problems of the form}
\begin{equation}\label{1.1}
	\min\limits_{x \in \mathbb{R}^n} F(x) = f(x) + r(x),
\end{equation}
{where $f$ is a smooth convex function and $r$ is a (possibly) non-smooth
	convex regularizer. A typical example is the expected risk minimization problem}
\[
f(x) := \mathbb{E}_{\xi \sim P}[\Lambda(x; \xi)],
\]
{where $\xi$ denotes a random data sample and $\Lambda(\cdot;\xi)$ is a convex loss function,
	or its finite-sum approximation}
\[
f_N(x) := \frac{1}{N}\sum_{i=1}^N \Lambda(x;\xi_i).
\]
{Such problems cover many standard tasks in machine learning and statistics, for example,  expected risk minimization problems such as logistic regression with $\ell_1$ or group-sparsity regularization. In these applications,
	$x$ represents the parameter vector to be learned from data and the regularization term  $r(x)$  is employed to induce structural properties, such as sparsity or group structure.}

In this setting, a natural and widely used approach is the stochastic gradient descent (SGD)
method~\cite{robbins1951stochastic}. In the finite-sum case $f(x) = \frac{1}{N}\sum_{i=1}^N f_i(x)$,
SGD draws randomly $i_k$ from $[N] := \{1,2,\dots,N\}$ and updates $x^{k+1}$ by
\[
x^{k+1} := x^k - \eta_k \nabla f_{i_k}(x^k),
\]
at each iteration. The advantage of the SGD method is that it only evaluates the gradient
of a single component function in each iteration, and hence the computational cost per
iteration is cheaper than that in the gradient descent method (GD). However, owing to the
variance unintentionally generated by random sampling, the SGD method converges slower than
the GD method. To overcome this drawback, various variance-reduction techniques have been
successively proposed; see
\cite{cutkosky2019momentum,dai2023variance,defazio2014saga,huang2021training,pham2020proxsarah,phan2024stochastic,traore2024variance,xiao2014proximal,xu2023momentum}.

Variance-reduction techniques inherit the advantage of low iteration cost of the SGD
method. Xiao et al.~\cite{xiao2014proximal} proposed the proximal stochastic
variance-reduced gradient (ProxSVRG) method which combines the SVRG~\cite{johnson2013accelerating}
variance-reduction technique with proximal mapping. The variance-reduction steps are as
follows:
\[
\left\{
\begin{aligned}
	&\textbf{Outer Loop:} \text{ For } s = 1, 2, \ldots: \\
	&\quad \tilde{x} = \tilde{x}_{s-1}, \quad \tilde{v} = \nabla F(\tilde{x}), \quad x_0 = \tilde{x}, \\
	&\textbf{Inner Loop:} \text{ For } k = 1, 2, \ldots, m: \\
	&\quad i_k \sim Q, \quad v_k = \frac{\nabla f_{i_k}(x_{k-1}) - \nabla f_{i_k}(\tilde{x})}{q_{i_k} n} + \tilde{v}.
\end{aligned}
\right.
\]

From the above, it can be seen that the ProxSVRG method requires computing an extra full
gradient $\nabla F$ every epoch. Further, Defazio et al.~\cite{defazio2014saga} proposed
the SAGA method which requires a full gradient in the first iteration and stores a history
of stochastic gradients in a matrix of size $N \times n$, where $N$ is the size of the
dataset and $n$ is the number of optimization variables.
{In large-scale regimes, however, computing a full gradient at every
	epoch or maintaining a gradient table of size $N\times n$ becomes prohibitively expensive
	in terms of both time and memory, especially when $N$ is very large.}
Thus, the ProxSVRG and SAGA methods are not suitable for large-scale data problems in general.

In order to overcome this difficulty, Dai et al.~\cite{dai2023variance} proposed the
S\mbox{-}PStorm algorithm, which employs a variance-reduction with momentum technique
\eqref{1.2} and a stabilized step-size strategy \eqref{1.3}. The algorithm is as follows:
\begin{empheq}[left=\empheqlbrace]{align}
	&\text{Sample } B_{k} = \{\xi_{k1}, \ldots, \xi_{km}\} \text{ (independent samples)}. \nonumber \\
	&\text{Compute } v_{k} = \frac{1}{m} \sum_{i=1}^m \nabla f(x_{k}; \xi_{ki}), \
	u_{k} = \frac{1}{m} \sum_{i=1}^m \nabla f(x_{k-1}; \xi_{ki}). \nonumber \\
	&\text{Update } d_k = v_{k} + (1 - \beta_k)(d_{k-1} - u_{k}). \tag{1.2} \label{1.2} \\
	&\text{Compute } y_k = \mathrm{prox}_{\alpha_k r}(x_k - \alpha_k d_k). \nonumber \\
	&\text{Update } x_{k+1} = x_k + \zeta\beta_k (y_k - x_k). \tag{1.3} \label{1.3} \nonumber
\end{empheq}
where $d_{k}$ is the gradient estimation, $\zeta \in (0,+ \infty)$, and $\beta_{k}$ is the
momentum coefficient. However, in the S\mbox{-}PStorm method, the step size $\alpha_k$ must
be fixed.

Variance-reduction algorithms mentioned above use fixed or diminishing step sizes; see
also \cite{duchi2011adaptive,huang2021training}. However, neither of these two approaches
can be efficient.

Recently, Tan et al.~\cite{tan2016barzilai} proposed the SVRG-BB algorithm that combines
the SVRG method with the Barzilai--Borwein (BB) step size~\cite{barzilai1988two}. The BB
step size uses past gradient information to adaptively calculate step sizes, avoiding
line search. The forms of BB step sizes are as follows:

BB1 step size (long step size)
\begin{equation}\label{longbb}
	\alpha_k^{\mathrm{BB1}} = \dfrac{\|s_k\|^2}{s_k^\top y_k},\tag{1.4}
\end{equation}
BB2 step size (short step size)
\begin{equation}
	\alpha_k^{\mathrm{BB2}} = \dfrac{s_k^\top y_k}{\|y_k\|^2}, \tag{1.5}\label{shortbb}
\end{equation}
where $s_k = x_k - x_{k-1}$ and $y_k = \nabla f(x_k) - \nabla f(x_{k-1})$. In practical
applications of the BB method, the step sizes $\alpha_k^{\mathrm{BB1}}$ and
$\alpha_k^{\mathrm{BB2}}$ are often used alternately.

The step size in \cite{tan2016barzilai} uses the BB1 step size~\eqref{longbb} as follows:
\[
\eta_k = \frac{1}{m} \cdot \| \tilde{x}_k - \tilde{x}_{k-1} \|_2^2 /
\big((\tilde{x}_k - \tilde{x}_{k-1})^\top (g_k - g_{k-1})\big),
\]
where $\eta_k$ is the step size, $g_k = \frac{1}{n} \sum_{i=1}^n \nabla f_i (\tilde{x}_k)$
and $m$ is the update frequency. In \cite{tan2016barzilai}, the numerical results show that
SVRG-BB is comparable to and sometimes even better than SVRG~\cite{johnson2013accelerating}
with best-tuned fixed step sizes, but as in \cite{dai2023variance,tan2016barzilai}, the
objective function $f(x)$ is required to be strongly convex.

{More recently, Yang et al.~\cite{Yang2023iSARAHBB} developed an inexact
	stochastic recursive gradient method with Barzilai--Borwein step sizes, called
	iSARAH-BB, for large-scale machine learning problems. They consider expected risk
	minimization and its finite-sum approximation, where each component function is convex and
	the objective is strongly convex. By incorporating a BB-type rule into the iSARAH
	framework, iSARAH-BB can adaptively compute dynamic step sizes and ensure convergence under these assumptions. Numerical experiments on
	$\ell_2$-regularized logistic regression show that iSARAH-BB is robust to the choice of the
	initial step size and competitive with several existing variance-reduced stochastic
	methods.}

{The variance-reduced stochastic methods mentioned above face several challenges in large-scale convex composite optimization: (1) Methods like ProxSVRG~\cite{xiao2014proximal} and SAGA~\cite{defazio2014saga} incur prohibitive per-epoch computation or memory costs; (2) momentum-based methods like S-PStorm avoid full gradients but rely on fixed step sizes; (3) adaptive BB-step methods like SVRG-BB \cite{tan2016barzilai} and iSARAH-BB~\cite{Yang2023iSARAHBB} enable dynamic step size, but they typically require the strong convexity of $F$ and do not handle general convex objectives with an explicit non-smooth regularizer $r(x)$, where $r(x)$ helps prevent overfitting.}

{To better address these challenges, motivated by \cite{dai2023variance,phan2024stochastic,tan2016barzilai,zhou2024adabb,Yang2023iSARAHBB}, we propose a stochastic proximal gradient method with adaptive step size and variance-reduction technique (PSGA) for efficiently solving the large-scale convex composite problem~\eqref{1.1} arising from regularized empirical risk minimization. In the PSGA method, we  incorporate a new variance-reduced stochastic gradient estimator that avoids computing full gradients every epoch and storing the history of gradient. We also employ a proximal operator which can handle the non-smooth regularization term $r(x)$. {A key component of the PSGA method is the adaptive step-size strategy based on the BB2 step size~\eqref{shortbb}. Since BB-type methods may diverge for general convex functions when the step size becomes too aggressive~\cite{zhou2024adabb}, we decrease the step size in the next iteration if it is too large and increase it if it is too small. This can avoid keeping the step size always small or always large, thereby ensuring fast convergence; see Step~5 in Algorithm~\ref{PSGA}. In addition, the objective function $F$ in our method is only  required to be convex.}
}

	Our contributions are summarized as follows.
	\begin{itemize}
		\item Different from  the assumption of strong convexity for the objective function in  \cite{dai2023variance,tan2016barzilai}, the objective function $f(x)$ for our method is only required  to be convex.
	\end{itemize}
	
	\begin{itemize}
		\item By adopting an adaptive step size strategy and  variance reduction technique, we avoid full gradient computations and historical gradient storage. In addition, the step size for our method is not necessarily fixed. At the same time, we prove that  the gradient estimation error converges to zero almost surely which implies the convergence in probability in \cite{dai2023variance}. This adaptive step size strategy also prevents the potential divergence of SVRG-BB\cite{tan2016barzilai} when applied to general convex functions.
	\end{itemize}
	
	\begin{itemize}
		\item Compared with the $O\left({\sqrt{\frac{\log k}{k}}}\right)$ convergence rate of the S-PStorm method in \cite{dai2023variance}, we achieve an improved rate of $O\left(\sqrt{\frac{1}{k}}\right)$ for our method.
	\end{itemize}
	
	\begin{itemize}
		\item We perform numerical experiments on Logistic regression and Lasso regression, demonstrating that our method achieves faster convergence and more accurate gradient estimation compared with  S-PStorm\cite{dai2023variance}, SAGA\cite{defazio2014saga}, RDA\cite{xiao2009dual}, Prox-SVRG\cite{xiao2014proximal}, and PStorm\cite{xu2023momentum} methods
	\end{itemize}
	
	The rest of this paper is organized as follows. In Section~\ref{Section2}, we introduce
	our algorithm. In Section~\ref{Section3}, we provide the definitions and assumptions
	required for the convergence proof and complete the proof. In Section~\ref{Section4}, we
	present our experimental results. Conclusions are given in Section~\ref{Section5}.

	\section{Algorithms}\label{Section2}
	{In this section, we present a proximal stochastic gradient algorithm (PSGA)
		for solving problem~\eqref{1.1}. The algorithm combines a variance-reduced stochastic
		gradient estimator with an adaptive Barzilai--Borwein (BB) step-size strategy and a
		proximal operator, so that it can efficiently handle large-scale problems with non-smooth
		regularization.} { In Step 3 of Algorithm~\ref{PSGA}, we implement a novel stochastic gradient estimator $\widetilde{\nabla}f(x_k)$. This design incorporates information from the current mini-batch and previous iterates, which can reduce the variance of the stochastic gradient without computing full-gradient at every epoch and storing the whole history of gradients.} {At the same time, we use an adaptive BB-type step-size rule in the PSGA method. This prevents the step size from staying small or large for a long time and thus contributes to faster practical convergence and improves the robustness; see Step 5 in Algorithm~\ref{PSGA}. } { We also apply the proximal operator to handle the non-smooth regularizer; see \eqref{2.1} in Algorithm ~\ref{PSGA}.}

	{To better illustrate the differences between our PSGA and other methods, we present a comparison in Table~\ref{tab:comparison_compact}, where $0 <\rho_{SVRG} <1.$}
	
	\begin{table}[htbp]
		\centering
		\small
		\caption{{Comparison of stochastic optimization methods}}
		\label{tab:comparison_compact}
		
		% tighten column spacing for this table only
		\begingroup
		\setlength{\tabcolsep}{3pt} % default is 6pt; reduce if too wide
		\renewcommand{\arraystretch}{1.05}
		
		\begin{tabular}{lccccc}
			\toprule
			& \rothead{\textbf{PSGA}\\ \textbf{(Proposed)}}
			& \rothead{\textbf{S-PStorm}\\ \textbf{\cite{dai2023variance}}}
			& \rothead{\textbf{SAGA}\\ \textbf{\cite{defazio2014saga}}}
			& \rothead{\textbf{Prox SVRG}\\ \textbf{\cite{xiao2014proximal}}}
			& \rothead{\textbf{RDA}\\ \textbf{\cite{xiao2009dual}}} \\
			\midrule
			\textbf{Step-size}
			& {Adaptive BB} & Fixed & Fixed & Fixed & Diminishing \\
			
			\textbf{Storage}
			& ${O}(n)$ & ${O}(n)$ & ${O}(Nn)$ & ${O}(n)$ & ${O}(n)$ \\
			
			\textbf{Full gradient}
			& {prob.\ $1/m$} & {Never} & Initial only & Per epoch & {Never} \\
			
			\textbf{Convexity}
			& Convex & Strong & Convex  & Strong & Strong \\
			
			\textbf{Rate}
			& \textbf{${O}(1/\sqrt{k})$}
			& ${O}(\sqrt{\log k/k})$
			& \shortstack{$O(1/k)$}
			& ${O}(\sqrt{\rho_{SVRG}^{k}})$
			& ${O}(\sqrt{\log k/k})$ \\
			
			\bottomrule
		\end{tabular}
		
		\endgroup
		
		\vspace{0.4em}
		\raggedright\footnotesize
		{}
	\end{table}
	
	\begin{remark}
		{In Table \ref{tab:comparison_compact}, N is the number of samples and n is the feature dimension.  ``Step-size'' denotes the step-size schedule, ``Storage'' represents the per-iteration memory complexity,  ``Full gradient'' describes when (or how often) a full gradient is computed, ``Convexity'' means the convexity of the objective function $F$ and  ``Rate" denotes the  convergence rate corresponding to the convexity assumption of $F$.}
	\end{remark}

	\newpage
	\begin{algorithm}[H]
		\caption{PSGA}
		\label{PSGA}
		\begin{algorithmic} % ??? [1]?????????
			\State \textbf{Step 1.} Choose initial point $x_0 = x_1 \in \mathbb{R}^n$, mini-batch size {$l \in \mathbb{N}^+$}, weight sequence $\{\theta_k\}_{k \geq 1} \in (0, 1)$ with $\theta_k = \frac{1}{k+1}$, step size sequence $\{\eta_k\} \in (0, +\infty)$ where $\eta_0 \geq \dfrac{1}{L}$, positive integer $m$, and $\delta_k = k$.
			
			\State \quad {Draw $l$ i.i.d.} samples $\{\xi_{k1}, \dots, \xi_{kl}\}$ from $\mathbb{P}$.
			
			\State \textbf{Step 2.} Compute
			\begin{equation*}
				{\mu_k = \frac{1}{l} \sum_{i=1}^l \nabla \Lambda(x_k; \xi_{ki}),}
			\end{equation*}
			\begin{equation*}
				{\nu_k = \frac{1}{l} \sum_{i=1}^l \nabla \Lambda(x_{k-1}; \xi_{ki}).}
			\end{equation*}
			
			\State \textbf{Step 3.} Compute
			\qquad $\widetilde{\nabla}f(x_k) = \mu_k$ \quad \text{if } $k=1$
			\qquad $\widetilde{\nabla}f(x_k) = \begin{cases}
				\nabla f(x_k) & \text{with prob. } 1/m, \\
				\mu_k + (1 - \theta_k)(\widetilde{\nabla}f(x_{k-1}) - \nu_k) & \text{with prob. } 1-1/m.
			\end{cases}$ \quad \text{if } $k > 1$
			
			\State \textbf{Step 4.} Compute
			\begin{equation}\label{case0}
				\tau_k= \dfrac{\langle \mu_k-\nu_k, x_k - x_{k-1} \rangle}{\|\mu_k - \nu_k\|^2}.
			\end{equation}
			
			\State \textbf{Step 5.} Set step size:
			\begin{equation}\label{case1}
				\text{If } \tau_k \geq \eta_{k-1}, \text{ set } \eta_k = \big(1 + \dfrac{1}{\tau_k}\big) \eta_{k-1},
			\end{equation}
			\begin{equation}\label{case2}
				\text{if } \eta_{k-1}/2 < \tau_k < \eta_{k-1}, \text{ set } \eta_k = \tau_k,
			\end{equation}
			\begin{equation}\label{case3}
				{\text{if } \tau_k \leq \eta_{k-1}/2, \text{ set } \eta_k = \dfrac{\eta_{k-1}}{2\sqrt{(\sqrt{m}+1)}}.}
			\end{equation}
			
			\State \textbf{Step 6.} Compute
			\begin{equation}\label{2.1}
				y_{k} = \mathrm{prox}_{\eta_k D(\cdot,x_k)} (x_k - \eta_k \widetilde{\nabla}f(x_k)),
			\end{equation}
			\begin{equation}\label{2.2}
				x_{k+1} = x_k + \delta_k \theta_{k} (y_k-x_k).
			\end{equation}
			
			\State \textbf{Step 7.} Update $k \gets k+1$ and return to Step 2.
		\end{algorithmic}
		
	\end{algorithm}

	\newpage
	\section{Convergence Analysis}\label{Section3}
	\noindent
	We begin this section by introducing some definitions, assumptions and lemmas.

	\begin{definition}\label{def3.1}
		(Surrogate function)\cite{phan2024stochastic}: A function \( D : \mathbb{R}^d \times \mathbb{R}^d \to \mathbb{R} \cup \{ +\infty \} \) is said to be a surrogate function of \( r : \mathbb{R}^d \to \mathbb{R} \cup \{ +\infty \} \) if
		
		\begin{enumerate}
			\item[(a)] \( D(y, y) = r(y) \) for all \( y \in \mathbb{R}^d \),
			\item[(b)] \( D(x, y) \geq r(x) \) for all \( x, y \in \mathbb{R}^d \).
		\end{enumerate}
	\end{definition}
	
	\begin{definition}\label{def3.2}
		(Almost surely)\cite{kolmogorov1956foundations}:
		An event $A$ is called almost surely (for short, a.s.) if $P(A) = 1$
	\end{definition}
	
	\begin{definition}\label{def3.3}
		\cite{ash2000probability}Let $\{A_n\}$ be the sequence of sets. The limit superior (or upper limit) is defined as
		\[
		\limsup A_n = \bigcap_{m=1}^\infty \bigcup_{n=m}^\infty A_n = \{\omega \mid \forall N, \exists n \geq N, \omega \in A_n\},
		\]
		and the limit inferior (or lower limit) is defined as
		\[
		\liminf A_n = \bigcup_{m=1}^\infty \bigcap_{n=m}^\infty A_n = \{\omega \mid \exists N, \forall n \geq N, \omega \in A_n\}.
		\]
		We set $\omega \in \limsup A_n$ as $\omega \in A_n$ infinitely often (abbreviated as $\omega \in A_n$ i.o.), meaning that $\omega$ belongs to $A_n$ for infinitely many $n$.
	\end{definition}
	
	\begin{lemma}\label{lem3.11}
		(Borel Cantelli Lemma)\cite{durrett2019probability}
		If the sum of the probabilities of the events $\{A_n\}$ is finite
		\[
		\sum_{n=1}^\infty P(A_n) < \infty,
		\]
		then the probability that infinitely many of them occur is 0, that is
		\[
		P\left(\limsup_{n \to \infty} A_n\right) = 0.
		\] (i.e., the probability that event $A_n$ occurs infinitely often is 0)
	\end{lemma}
	
	\begin{lemma}\label{lem3.12}
		(Markov's inequality)\cite{lin2010probability}
		If $\varphi$ is a non-decreasing non-negative function, $X$ is a (not necessarily nonnegative) random variable, and $\varphi(a) > 0$, then
		
		\[
		\mathrm{P}(X \geq a) \leq \frac{\mathrm{E}(\varphi(X))}{\varphi(a)}.
		\]
	\end{lemma}
	
	For Problem (\ref{1.1}), the following assumptions are required:
	
	\begin{assumption}\label{ass3.4}
		$f$ is convex over $\mathbb{R}^n$ and $r$ is convex and closed over $\mathbb{R}^n$.
	\end{assumption}
	
	\begin{assumption}\label{ass3.5}
		There exists a constant $L > 0$ such that, for any $(x, y) \in \mathbb{R}^n \times \mathbb{R}^n$ and any $\xi \sim \mathbb{P}$, it holds that
		\begin{center}
			$\|\nabla \Lambda(x, \xi) - \nabla \Lambda(y, \xi)\| \leq $L$\|x - y\|$,
		\end{center}
		i.e., $f(x)$ is $L$-smooth.
		
	\end{assumption}
	
	\begin{assumption}
		There exists $G_r > 0$ such that, for all $k \geq 1$,
		\[
		\mathbb{P}\{\|g_r\|_2 \leq G_r, \;  g_r \in \partial r(x_k)\} = 1.
		\]
	\end{assumption}
	
	To ensure the convergence of Algorithm \ref{PSGA}, we require the following assumption:
	\begin{assumption}\label{ass3.6}
		\item[(a)]For all $k \geq 1$, $\mathbb{E}_{\xi \sim \mathcal{P}}[\nabla \Lambda(x_k; \xi) \mid \mathcal{F}_k] = \nabla f(x_k)$,
		where $\mathbb{E}_{\xi \sim \mathcal{P}}[\nabla \Lambda(x_k; \xi) \mid \mathcal{F}_k]$ denotes that the expected value of the stochastic gradient $\nabla \Lambda(x_k; \xi)$ over the sample distribution $P$, conditioned on the historical information $\mathcal{F}_k.$
		\item[(b)] There exists $\sigma > 0$ such that, for all $k \geq 1$,
		\[
		\mathbb{P}_{\xi \sim \mathcal{P}}\{\|\nabla \Lambda(x_k, \xi) - \nabla f(x_k)\| \leq \sigma \mid \mathcal{F}_k\} = 1.
		\]
	\end{assumption}
	
	Similar to  problem (\ref{1.1}), for Surrogate function $D(x, y)$, we need the following assumption.
	
	\begin{assumption}\label{ass3.7}
		\item[(a)] For every \( x \), \( D(x, \cdot) \) is continuous in \( y \).
		\item[(b)] For every \( y \), \( D(\cdot, y) \) is lower semicontinuous and convex.
		\item[(c)] There exists a function \( {c} : \mathbb{R}^d \times \mathbb{R}^d \to \mathbb{R} \ \) such that for every \( y \in \mathbb{R}^d \), \( {c}(\cdot, y) \) is continuously differentiable at \( y \) with \( \nabla{c}(\cdot, y) (y) = 0 \), and the approximation error satisfies
		\[
		D(\cdot, y) - r(\cdot) \leq {c}(\cdot, y).
		\]
	\end{assumption}

	\subsection{Convergence Analysis}
	\noindent
	
	In order to prove the convergence of Algorithm \ref{PSGA}, we need the following lemmas:
	
	\begin{lemma}\label{lem3.8}
		\cite{robbins1971convergence}
		Let $\{Y_k\}$, $\{Z_k\}$, and $\{W_k\}$ be three sequences of random variables and let $\mathcal{F}_k$ be sets of random variables such that $\mathcal{F}_k \subset \mathcal{F}_{k+1}$ for all $k$. Assume that
		
		\begin{enumerate}
			\item[(a)] The random variables $\{Y_k\}$, $\{Z_k\}$, and $\{W_k\}$ are nonnegative and are functions of random variables in $\mathcal{F}_k$;
			\item[(b)] $\mathbb{E}[Y_{k+1} \mid \mathcal{F}_k] \leq Y_k - Z_k + W_k$ for each $k$;
			\item[(c)] $\sum_{k=0}^\infty W_k < +\infty$ with probability 1.
		\end{enumerate}
		
		Then, $\sum_{k=0}^\infty Z_k < +\infty$, and $\{Y_k\}$ converges to a nonnegative random variable, almost surely.
	\end{lemma}
	
	The following two lemmas play a important role in the convergence analysis of Algorithm \ref{PSGA}.
	
	\begin{lemma}\label{lem3.9}
		Suppose that Assumptions \ref{ass3.4}-\ref{ass3.6} hold. Let $\{x_k\}$ be the sequence generated by Algorithm \ref{PSGA} and $\Psi_k:=\|\widetilde{\nabla}f(x_{k-1}) - \nabla f(x_{k-1})\|^2$, and let $\mathbb{E}_k$ denote the conditional expectation on $\mathcal{F}_k$.  Then
		\begin{align}\label{in3.9}
			\mathbb{E}_k \left\|\widetilde{\nabla}f(x_k) - \nabla f(x_k) \right\|^2 \leq \Psi_k + 4{L^2} \|x_k - x_{k-1}\|^2 + 2{\theta_k}^2\sigma^2.
		\end{align}
	\end{lemma}
	
	\noindent  Proof.
	From the Algorithm $\ref{PSGA}$ and the definition of $\widetilde{\nabla}f(x_k)$, we get
	\begin{align}
		\mathbb{E}_k \|\widetilde{\nabla}f(x_k) &- \nabla f(x_k)\|^2 \notag \\
		&= \left(1-\frac{1}{m}\right)\mathbb{E}_k\|\mu_k + (1-\theta_k)(\widetilde{\nabla}f(x_{k-1})-\nu_k) - \nabla f(x_k)\|^2 \notag \\
		&\leq \left(1-\frac{1}{m}\right)\Bigl(\mathbb{E}_k\|\mu_k-\nabla f(x_k) + (1-\theta_k)(\nabla f(x_{k-1}) - \nu_k)\|^2 \notag \\
		&\quad + \mathbb{E}_k \|\widetilde{\nabla}f(x_{k-1}) - \nabla f(x_{k-1})\|^2\Bigr) \notag \\
		&\leq \left(1-\frac{1}{m}\right)\Bigl(\mathbb{E}_k\bigl[2\|\mu_k - \nu_k + \nabla f(x_{k-1}) - \nabla f(x_k)\|^2 \notag \\
		&\quad + 2\theta_k^2\|\nu_k - \nabla f(x_{k-1})\|^2\bigr]\Bigr) + \left(\frac{m-1}{m}\right) \Psi_k \label{3.0} \\
		&\leq \left(1-\frac{1}{m}\right)\Bigl[4L^2\|x_k - x_{k-1}\|^2 + \Psi_k + 2\theta_k^2\sigma^2\Bigr] \label{3.1}  \\
		&\leq \Psi_k + 4L^2\|x_k - x_{k-1}\|^2 + 2\theta_k^2\sigma^2, \label{3.1.1}
	\end{align}
	where the second inequality follows from the inequality $\|a+b\|^2\leq2\|a\|^2+2\|b\|^2$, the third inequality is obtained from Assumptions $\ref{ass3.5}$ and $\ref{ass3.6}(b)$, and the fourth inequality is due to $0 < 1-\dfrac{1}{m} < 1$. This completes the proof.
	
	\begin{lemma}\label{the3.10}
		Let $\{x_k\}$ be the sequence generated by  Algorithm \ref{PSGA}. Suppose that Assumptions \ref{ass3.4}-\ref{ass3.6} hold and
		\begin{equation}\label{new3.3}
			\eta_{k} \leq \frac{k+1}{4(\sqrt{m}+1)\delta_k L}.
		\end{equation}
		Then
		\item[(a)] The sequence $\{\|x_k - x_{k-1}\|^2\}$ has a finite sum almost surely. \label{the3.10(b)}
		\item[(b)] The sequence $\{F(x_k)\}$ converges almost surely. \label{the3.10(a)}
	\end{lemma}
	
	\noindent  Proof.
	(a) Since $r$ is convex and $\delta_k \theta_{k} < 1$, from (\ref{2.2}) we have
	\begin{equation}\label{3.2}
		r(x_{k+1}) \leq \delta_k \theta_{k}r(y_k)+(1-\delta_k \theta_{k})r(x_k)
		\leq \delta_k \theta_kD(y_{k},x_k)+(1-\delta_k \theta_k)r(x_k).
	\end{equation}
	The $L$-smoothness of $f$ yields
	\begin{equation}\label{3.3}
		f(x_{k+1}) \leq f(x_k) + \langle \nabla f(x_k), x_{k+1} - x_k \rangle + \frac{L}{2} \|x_{k+1} - x_k\|^2.
	\end{equation}
	In view of the definition of the proximal operator, we obtain
	\begin{align}\label{new3.6}
		D(y_k, x_k) + \frac{1}{2\eta_{k}}\|y_{k} - (x_k-\eta_{k}\widetilde{\nabla}f(x_k))\|^2 \leq D(x_k,x_k) + \frac{1}{2\eta_{k}}\|\eta_{k}\widetilde{\nabla}f(x_k)\|^2.
	\end{align}
	Note that
	\begin{align}\label{new3.7}
		\|y_{k} - (x_k-\eta_{k}\widetilde{\nabla}f(x_k))\|^2 = \|y_{k} - x_k\|^2 + 2\langle \eta_{k}\widetilde{\nabla}f(x_k),y_{k}-x_k \rangle + \| \eta_{k}\widetilde{\nabla}f(x_k)\|^2.
	\end{align}
	Combining (\ref{new3.7}) with (\ref{new3.6}), we have
	\begin{equation}\label{3.4}
		\langle \widetilde{\nabla}f(x_k), y_{k} - x_k \rangle + D(y_{k}, x_k) + \frac{1}{2\eta_k} \|y_{k} - x_k\|^2 \leq D(x_k, x_k) = r(x_k).
	\end{equation}
	Multiplying both sides by $\delta_k \theta_{k}$ in (\ref{3.4}), we get
	
	\begin{align}\label{3.5}
		\langle \widetilde{\nabla}f(x_k),x_{k+1}-x_k \rangle +\delta_k \theta_k D(y_k,x_k) + \frac{1}{2\delta_k \theta_k \eta_{k}} \|x_{k+1}-x_k\|^2 \leq \delta_k \theta_k r(x_k).
	\end{align}
	
	From (\ref{3.2}), (\ref{3.3}) and (\ref{3.5}), we deduce that
	\begin{align}
		F(x_{k+1}) + (\frac{1}{2\delta_k \theta_k\eta_{k}}+\frac{L}{2})\|x_{k+1}-x_k\|^2
		&\leq F(x_{k}) + \langle \nabla f(x_k) - \widetilde{\nabla}f(x_k) ,x_{k+1}-x_k \rangle \notag\\
		&\leq F(x_{k}) + \frac{\xi}{2}\|\widetilde{\nabla}f(x_k)-\nabla f(x_k)\|^2 + \frac{2}{\xi}\|x_{k+1} - x_k\|^2,\label{3.7}
	\end{align}
	where $\xi>0$ is arbitrary.
	By taking expectation in (\ref{3.7}), conditioned on $\mathcal{F}_k$, we have
	\begin{align}\label{3.8}
		\mathbb{E}_k[F(x_{k+1}) + (\frac{1}{2\delta_k \theta_k\eta_{k}}-\frac{L}{2}-\frac{2}{\xi})\|x_{k+1}-x_k\|^2] \leq F(x_k)+ \frac{\xi}{2} \mathbb{E}_k \|\widetilde{\nabla}f(x_k)-\nabla f(x_k)\|^2.
	\end{align}
	
	Letting $J = \frac{1}{m}, \Pi_{\Psi} = \frac{(m-1)L^2}{m}$ and $ \Pi = L^2$ in (\ref{3.8}), using Lemma \ref{lem3.9} we get
	\begin{align}{2}\label{3.8a}
		&\mathbb{E}_k \bigg[ F(x_{k+1})
		+ \left( \frac{1}{2\delta_k\theta_k\eta_k} - \frac{L}{2} - \frac{2}{\xi} \right) \|x_{k+1} - x_k\|^2
		+ \frac{\xi}{2J} \Psi_{k+1} \bigg] \nonumber\\
		&\quad \leq F(x_k) + \frac{\xi}{2J} \Psi_k
		+ \left( \frac{\xi}{J} \Pi_{\Psi} + J \Pi \right) \|x_{k+1} - x_k\|^2
		+ \left( 1 + \frac{1}{J} \right) \xi \theta_k^2 \sigma^2.
	\end{align}
	
	Setting
	\begin{align}\label{3.8b}
		\gamma_{k+1} = F(x_{k+1})+(\dfrac{1}{2\delta_k\theta_k\eta_{k}}+\dfrac{L}{2}-\dfrac{2}{\xi})\|x_{k+1}-x_k\|^2+\dfrac{\xi}{2J}\Psi_{k+1}, \end{align}
	and from (\ref{3.8a}) we obtain

	\begin{align}
		\mathbb{E}_k \gamma_{k+1}
		&\leq \gamma_{k} -(\frac{1}{2\delta_k\theta_k\eta_{k}}-\frac{L}{2}-\frac{2}{\xi}-\frac{\xi \Pi}{2}-\frac{\xi \Pi_{\Psi}}{2J})\|x_k-x_{k-1}\|^2 + (1+\frac{1}{J})\xi\theta_k^2\sigma^2 \notag\\
		&\leq \gamma_{k} -(2(\sqrt{m}+1)L-\frac{L}{2}-\frac{2}{\xi}-\frac{\xi \Pi}{2}-\frac{\xi \Pi_{\Psi}}{2J})\|x_k-x_{k-1}\|^2 + (1+\frac{1}{J})\xi\theta_k^2\sigma^2 \label{3.10},
	\end{align}
	
	where (\ref{3.10}) is obtained by the inequality (\ref{new3.3}). Setting $\xi = \dfrac{2}{\sqrt{m}L}$ we have
	\begin{align*}
		2(\sqrt{m}+1)L-\frac{L}{2}-\frac{2}{\xi}-\frac{\xi \Pi}{2}-\frac{\xi \Pi_{\Psi}}{2J}
		= 2(\sqrt{m}+1)L - \frac{L}{2} - 2\sqrt{m}L = \frac{3}{2}L > 0
	\end{align*}

	Since $\sum_{k=1}^{+\infty}\theta_{k}^2 = \sum_{k=1}^{+\infty}\dfrac{1}{(k+1)^2} < +\infty$, from Lemma \ref{lem3.8} we obtain
	
	$$\sum_{k=1}^{+\infty}\|x_k-x_{k-1}\|^2 < +\infty \quad \text{a.s.}$$ and \{$\gamma_k$\} converges to a non-negative random variable $\gamma_{\infty}$ almost surely.
	
	(b) Combining Lemma \ref{lem3.8} and Lemma \ref{lem3.9}, we  get that $\Psi_k$ has a finite sum almost surely. Thus, from (\ref{3.8b}) we deduce that \{$F(x_k)$\} converges to $\gamma_{\infty}$ almost surely.
	
	In the following, we present the error between stochastic gradient estimation $\widetilde{\nabla}f(x_k)$ and the true gradient $\nabla f(x_k)$.

	\begin{theorem}\label{the3.13}
		Suppose Assumptions \ref{ass3.4}-\ref{ass3.6} hold. Let $\{x_k\}$ be the sequence generated by Algorithm \ref{PSGA} and $\Psi_{k+1}:=\|\widetilde{\nabla}f(x_{k}) - \nabla f(x_{k})\|^2$. Then
		\begin{align}\notag
			\lim_{k \to +\infty}  [\widetilde{\nabla}f(x_k)-\nabla f(x_k)] = 0\quad a.s.
		\end{align}
	\end{theorem}
	
	\noindent  Proof.
	By taking the total expectation in (\ref{3.10}), we obtain
	\begin{equation}\label{3.12}
		\mathbb{E} \gamma_{k+1} \leq \mathbb{E} \gamma_{k} - G\mathbb{E}\|x_k-x_{k-1}\|^2 + (1+\frac{1}{J})\xi\theta_k^2\sigma^2,
	\end{equation}
	where $G = 2(\sqrt{m}+1)L-\frac{L}{2}-\frac{2}{\xi}-\frac{\xi \Pi}{2}-\frac{\xi \Pi_{\Psi}}{2J}$.
	For any $K\geq 1$, we sum the inequality (\ref{3.12}) for $k=1,2,\cdots,K$ to obtain
	\begin{equation}\label{new3.17}
		\sum_{k=1}^{K}G\mathbb{E}\|x_k-x_{k-1}\|^2 \leq \mathbb{E}\gamma_1 - F_* +(1+\frac{1}{J})\sum_{k=1}^{K}\xi\theta_k^2\sigma^2,
	\end{equation}
	by using the fact that $ F_* \leq F_{K+1} \leq \gamma_{K+1}$.
	Since $(1+\frac{1}{J})\sum_{k=1}^{K}\xi\sigma^2\theta_k^2 \leq G_0:=(1+\frac{1}{J})\xi\sigma^2\frac{\pi^2}{6}$, from (\ref{new3.17}) we have
	\begin{equation}\label{3.14}
		\sum_{k=1}^{K}\mathbb{E}\|x_k-x_{k-1}\|^2 \leq \frac{\mathbb{E}\gamma_1 - F_* + G_0}{G},
	\end{equation}
	which implies that \{$\mathbb{E}\|x_k-x_{k-1}\|^2$\} has a finite sum.
	
	Next let us think of the \{$\mathbb{E}\|\widetilde{\nabla}f(x_k)-\nabla f(x_k)\|^2$\}.
	By taking the total expectation in (\ref{3.1}), we have
	\begin{equation}\label{new3.20}
		\mathbb{E} \Psi_k \leq 4(\frac{1}{J}-1)L^2\mathbb{E}\|x_k-x_{k-1}\|^2 + \frac{\mathbb{E}\Psi_k-\mathbb{E} \Psi_{k+1}}{J} + 2(\frac{1}{J}-1) \theta_k^2\sigma^2.
	\end{equation}
	Similarly, taking the total expectation in (\ref{in3.9}) to obtain
	\begin{equation}\label{3.15}
		\mathbb{E}\|\widetilde{\nabla}f(x_k)-\nabla f(x_k)\|^2 \leq \mathbb{E}\Psi_k + 4{L^2} \mathbb{E}\|x_k - x_{k-1}\|^2 + 2{\theta_k}^2\sigma^2.
	\end{equation}
	
	Combining  (\ref{3.15}) with (\ref{new3.20}), we get
	
	\begin{equation}\label{3.22}
		\mathbb{E}\|\widetilde{\nabla}f(x_k)-\nabla f(x_k)\|^2 \leq \frac{4L^2}{J} \mathbb{E}\|x_k - x_{k-1}\|^2 + \frac{\mathbb{E}\Psi_k-\mathbb{E} \Psi_{k+1}}{J}+ \frac{2}{J}\theta_k^2\sigma^2.
	\end{equation}
	
	Summing the inequality (\ref{3.22}) for $k=1,2,\cdots,K$, we have
	\begin{align}
		\sum_{k=1}^{K}\mathbb{E}\|\widetilde{\nabla}f(x_k)-\nabla f(x_k)\|^2
		&\leq \frac{\mathbb{E}\Psi_0-\mathbb{E}\Psi_{K+1}}{J}+\frac{4L^2}{J}\sum_{k=0}^{K}\mathbb{E}\|x_k-x_{k-1}\| +\frac{2}{J}\sum_{k=0}^{K}\theta_{k}^2\sigma^2 \notag \\
		&\leq \frac{4L^2}{J}\sum_{k=0}^{K}\mathbb{E}\|x_k-x_{k-1}\|+\frac{\pi^2}{3J} \label{eq:3.18},
	\end{align}
	where the second inequality follows from $\mathbb{E}\Psi_0=0$, $\mathbb{E}\Psi_K \geq 0$, and $\sum_{k=0}^{K}\theta_{k}^2 \leq \sum_{k=0}^{+\infty}\theta_{k}^2 = \frac{\pi^2}{6}$. Therefore,
	\begin{align}\label{eq3.24}
		\sum_{k=0}^{+\infty}\mathbb{E}\|\widetilde{\nabla}f(x_k)-\nabla f(x_k)\|^2<+\infty.
	\end{align}
	
	Set \( Y_k = \|\widetilde{\nabla}f(x_k) - \nabla f(x_k)\| \).
	% we have
	%		\[\sum_{k=0}^{+\infty} \mathbb{E}[Y_k^2] < +\infty.\]
	%	To prove \( Y_k \to 0 \) a.s. , note that by \cite{durrett2019probability}, \( X_n \to 0 \) a.s. if and only if for all \( \epsilon > 0 \),
	%	\[ P(|X_n| > \epsilon \text{ i.o.}) = 0. \]
	%		So we need to prove
	%		$$ P\left(\limsup_{k\to+\infty} \{Y_k > \epsilon\}\right) = 0 \quad \text{for all } \epsilon > 0. $$
	For any given $\epsilon > 0$, define the events
	$$ A_k = \{Y_k > \epsilon\}. $$
	In view of Lemma \ref{lem3.12}, we have
	$$ P(A_k) = P(Y_k > \epsilon) \leq \frac{\mathbb{E}[Y_k^2]}{\epsilon^2}. $$
	{Thus, using (\ref{eq3.24}) we get}
	$$ \sum_{k=0}^{+\infty} P(A_k) \leq \frac{1}{\epsilon^2}\sum_{k=0}^{+\infty} \mathbb{E}[Y_k^2] < +\infty. $$
	
	Thus, from Lemma \ref{lem3.11} we obtain
	$$ P\left(\limsup_{k\to+\infty} A_k\right) = 0, $$
	which means
	$$ P(Y_k > \epsilon \text{ i.o}) = 0. $$
	Since $ Y_k \to 0 $ a.s. if and only if for all $ \epsilon > 0 $, $ P(|Y_k| > \epsilon \text{ i.o.}) = 0. $(see for example \cite{durrett2019probability}), $ Y_k \to 0 $ a.s. and hence
	%Since $\epsilon > 0$ was arbitrary, we conclude
	%$$$ P\left(\limsup_{k\to+\infty} \{Y_k > \epsilon\}\right) = 0 \quad \text{for all } \epsilon > 0. $$
	%i.e.
	%$$ P\left(\lim_{k\to+\infty} Y_k = 0\right) = 1. $$

	\begin{equation}\notag
		\lim_{k \to +\infty}  [\widetilde{\nabla}f(x_k)-\nabla f(x_k)] = 0, \quad \text{a.s.}
	\end{equation} which completes the proof.
	
	The following theorem  establishes the variance reduction property of the stochastic gradient estimator.
	
	\begin{theorem}\label{rem3.14}
		Suppose Assumptions \ref{ass3.4}-\ref{ass3.6} hold. Let $\{x_k\}$ be the sequence generated by Algorithm \ref{PSGA}. Then
		\begin{equation}\notag
			\min_{k=1,2,...K}\mathbb{E}\|\widetilde{\nabla}f(x_k)-\nabla f(x_k)\|^2 \leq \frac{G_2}{K}.
		\end{equation}
	\end{theorem}
	
	\noindent  Proof.
	
	Combining (\ref{new3.17}) with (\ref{eq:3.18}), we have
	\begin{equation}\label{in3.30}
		\sum_{k=1}^{K}\mathbb{E}\|\widetilde{\nabla}f(x_k)-\nabla f(x_k)\|^2 \leq \frac{4L^2(\mathbb{E}\gamma_1 - F_* + G_0)}{GJ}+\frac{\pi^2}{3J}.
	\end{equation}
	Setting $G_2 = \dfrac{4L^2(\mathbb{E}\gamma_1 - F_* + G_0)}{GJ}+\dfrac{\pi^2}{3J}$ , we obtain
	
	\begin{equation*}
		\min_{k=1,2,...K}\mathbb{E}\|\widetilde{\nabla}f(x_k)-\nabla f(x_k)\|^2 \leq \frac{G_2}{K},
	\end{equation*}
	which completes the proof.

	Next, we present the convergence result of Algorithm \ref{PSGA}.
	\begin{theorem}\label{the3.15}
		Let $\{x_k\}$ be the sequence generated by  Algorithm \ref{PSGA}. Suppose Assumptions \ref{ass3.4}-\ref{ass3.7} hold, and
		\begin{equation}\notag
			\eta_{k} \leq \frac{k+1}{4(\sqrt{m}+1)\delta_k L},
		\end{equation}
		then the limit point of $\{x_k\}$ is an optimal point of F almost surely.
	\end{theorem}
	
	\noindent  Proof.
	In view of Lemma \ref{the3.10(b)} (a) and Theorem \ref{the3.13}, we obtain
	\begin{equation}\notag
		\lim_{k \to +\infty}  [\widetilde{\nabla}f(x_k)-\nabla f(x_k)] = 0 \quad \text{a.s.} \quad
		and \lim_{k \to +\infty} [ x_k - x_{k-1}] = 0 \quad \text{a.s.}
	\end{equation}
	Let $x^*$ be a limit point of \{$x_k$\}. Then there exists a subsequence \{$x_{k_i}$\} of \{$x_k$\} such that $x_{k_i}\to x^*$($i \to +\infty$). In view of (\ref{2.1}), we have
	\begin{align}\label{3.28}
		0 \in \frac{1}{\eta_{k_i}}(y_{k_i}-x_{k_i}+\eta_{k_i}\widetilde{\nabla}f(x_{k_i})) + \partial D(\cdot,x_{k_i})(y_{k_i}).
	\end{align}
	Using the definition of $\partial D(\cdot,x_{k_i})(y_{k_i})$, from (\ref{3.28}) we get
	\begin{align}\label{3.29}
		D(x,x_{k_i})-D(y_{k_i},x_{k_i}) \geq \langle -\frac{1}{\eta_{k_i}}(y_{k_i}-x_{k_i}+\eta_{k_i}\widetilde{\nabla}f(x_{k_i})), x-x_{k_i} \rangle , \forall x \in R ^{n}.
	\end{align}
	From (\ref{2.2}), we obtain
	\begin{align}\notag
		\|y_k-x_k\|^2=\|\frac{x_{k+1}-x_k}{\delta_k \theta_k}\|^2=\frac{(k+1)^2}{k^2}\|x_{k+1}-x_k\|^2,
	\end{align}
	which  by Theorem \ref{the3.10}(a) implies
	\begin{align}\notag
		\lim_{k \to +\infty}\|y_k-x_k\|^2 = 0 \quad a.s.
	\end{align}
	
	Letting $x=x^*$ in (\ref{3.29}) and then taking superior limit yields
	\begin{align}\label{3.34}
		\limsup_{i \to +\infty} D(y_{k_i}, x_{k_i})\leq r(x^*),
	\end{align}
	being $D(x,\cdot)$  continuous.
	In view of the lower semicontinuity of $D(\cdot,y)$, from (\ref{3.34}) we get
	\begin{align}\notag
		\lim_{i \to +\infty} D(y_{k_i}, x_{k_i})=r(x^*).
	\end{align}
	Now letting $i \rightarrow +\infty$ in (\ref{3.29}), and hence we get
	\begin{align}\label{3.36}
		r(x^*)\leq - \langle \nabla f(x^*),x-x^* \rangle+D(x,x^*).
	\end{align}
	Since $f$ is $L-$smooth,
	\begin{align}\label{3.37}
		f(x^*) \leq f(x) - \langle \nabla f(x^*), x - x^* \rangle + \frac{L}{2} \|x - x^*\|^2.
	\end{align}
	Combining (\ref{3.36}) and (\ref{3.37}), from Assumption \ref{ass3.7}(c) we get
	\begin{align}\notag
		F(x^*) \leq F(x)+D(x,x^*)-r(x)+\frac{L}{2}\|x-x^*\|^2
		\leq  F(x)+c(x,x^*)+\frac{L}{2}\|x-x^*\|^2.
	\end{align}
	Therefore, $x^*$ is the minimizer of
	\begin{center}
		$\min\limits_{x \in \mathbb{R}^n} F(x)+c(x,x^*) + \frac{L}{2}\|x-x^*\|^2.$
	\end{center}
	Thus,
	\begin{center}
		$0 \in \partial F(x^*) + \nabla c(\cdot,x^*)(x^*) = \partial F(x^*),$
	\end{center}
	where $\nabla c(\cdot,x^*)(x^*)=0$ is due to Assumption \ref{ass3.7}(c). As a result, $x^*$ is the optimal point of $F$ almost surely.

	The following  lemma plays an important role in  proving the convergence rate for our method.
	
	\begin{lemma} \label{lem3.18}
		{Let $\{\eta_k\}$ be the sequence generated by  Algorithm \ref{PSGA}. Suppose Assumptions \ref{ass3.5} and \ref{ass3.6} hold, then
		\begin{align}\label {in3.20}
			\eta_k \geq C_0 := \dfrac{1}{4(\sqrt{m}+1)L},\,\forall\,k\geq 0.
		\end{align}}
		
	\end{lemma}
	
	\noindent  {Proof. The proof will be divided into three steps.}
	
	\textbf{Step 1.} $\tau_k \geq \dfrac{1}{L}.$
	
	Using Assumption \ref{ass3.6}, from \cite[Theorem 2.1.5]{Nesterov} we obtain
	\begin{align}\label{3.45}
		\langle\nabla \Lambda(x,\xi) - \nabla \Lambda(y,\xi),x-y\rangle \geq \frac{1}{L}\| \nabla \Lambda(x,\xi) - \nabla \Lambda(y,\xi) \|^2,\,\forall\,(x, y) \in \mathbb{R}^n.
	\end{align}
	In view of Step 4 of Algorithm \ref{PSGA}, we have
	
	\begin{align}
		\tau_{k}&= \dfrac{\langle \mu_k-\nu_k, x_k - x_{k-1} \rangle}{\|\mu_k - \nu_k\|^2}\notag \\
		&= \dfrac{l\sum_{i=1}^{l}\langle \nabla \Lambda(x_k,\xi_{ki}) - \nabla \Lambda(x_{k-1},\xi_{ki}),x_k-x_{k-1}\rangle}{\|\sum_{i=1}^{l}(\nabla \Lambda(x_k,\xi_{ki}) - \nabla \Lambda(x_{k-1},\xi_{ki}))\|^2} \notag \\
		&\geq \dfrac{l\sum_{i=1}^{l}\|\nabla \Lambda(x_k,\xi_{ki}) - \nabla \Lambda(x_{k-1},\xi_{ki})\|^2}{L\|\sum_{i=1}^{l}(\nabla \Lambda(x_k,\xi_{ki}) - \nabla \Lambda(x_{k-1},\xi_{ki}))\|^2} \notag \\
		&\geq \frac{1}{L},
	\end{align}
	where the first inequality follows from (\ref{3.45}) and the second one is due to {$l\sum_{i=1}^{l}\|a_i\|^2 \geq \|\sum_{i=1}^{l}a_i\|^2.$}
	
	\textbf{Step 2.} {If $ \eta_i \geq \dfrac{1}{2\sqrt{(\sqrt{m}+1)L}}$ for some $i$, then $\eta_k \geq \dfrac{1}{4{(\sqrt{m}+1)}L}$ for any $k\geq i. $}
	
	{Using the proof by induction, we only need to prove that the conclusion holds when $k=i+1$. According to Step 4 of Algorithm \ref{PSGA}, we will take into account three different situations. If $\tau_{i+1} \geq \eta_i$ , then  $\eta_{i+1}\geq\eta_{i}\geq \dfrac{1}{2\sqrt{(\sqrt{m}+1)}L}\geq\dfrac{1}{4(\sqrt{m}+1)L}$. If $\eta_{i}/2 < \tau_{i+1} < \eta_{i}$, then $\eta_{i+1} = \tau_{i+1} \geq \dfrac{1}{L}\geq\dfrac{1}{4(\sqrt{m}+1)L}$ by using  Step 1.  If $\tau_{i+1} \leq \eta_{i}/2$, then $\eta_{i+1} = \dfrac{\eta_{i}}{2\sqrt{(\sqrt{m}+1)}} \geq \dfrac{1}{4(\sqrt{m}+1)L}$.}
	
	\textbf{Step 3.}{ $\eta_k \geq C_0 := \dfrac{1}{4(\sqrt{m}+1)L},\,\forall\,k\geq 0.$}
	
	{Let $j\geq1$ be the smallest integer such that $\tau_j < \eta_{j-1},$ which means that  $\tau_k < \eta_{k-1}$ for any $k\geq j$.}
	
	{If $j=1$, which means $\tau_1 < \eta_0,$  then from (\ref{case2}) and (\ref{case3}) we have $\eta_1 = \tau_1 \geq \dfrac{1}{L}\geq\dfrac{1}{2\sqrt{(\sqrt{m}+1)}L}$ or $\eta_1= \dfrac{\eta_0}{2\sqrt{(\sqrt{m}+1)}} \geq \dfrac{1}{2\sqrt{(\sqrt{m}+1)}L}$.  Therefore, using Step 2 we deduce that $\eta_k \geq \dfrac{1}{4(\sqrt{m}+2)L}$ for any $k\geq1$.}
	
	{Suppose now that $j>1$. For $1\leq k\leq j-1$, $\tau_k \geq \eta_{k-1}$  and hence from (\ref{case1}) we obtain $\eta_k > \eta_{k-1} \geq \eta_0 \geq \dfrac{1}{L} \geq \dfrac{1}{4(\sqrt{m}+1)L}$. For  $ k= j$, $ \tau_{j}<\eta_{j-1}$ and hence from(\ref{case2}) and (\ref{case3}) we have $\eta_j= \tau_j \geq \dfrac{1}{L}\geq\dfrac{1}{2\sqrt{(\sqrt{m}+1)}L}$ or $\eta_j=\dfrac{\eta_{j-1}}{2\sqrt{(\sqrt{m}+1)}} \geq\dfrac{\eta_0}{2\sqrt{(\sqrt{m}+1)}} \geq \dfrac{1}{2\sqrt{(\sqrt{m}+1)}L}.$ Thus, from Step 2 we obtain that $\eta_k \geq \dfrac{1}{4(\sqrt{m}+1)L}$ for $k\geq j$. As a result, $\eta_k \geq \dfrac{1}{4(\sqrt{m}+1)L}$ for any $k\geq1$.}
	
	In order to achieve the convergence rate, we need the following additional assumption (\cite{Phan Duy Nhat}).
	\begin{assumption}\label{ass3.11}
		For any bounded subset $\Omega$ of $\mathbb{R}^d$, there exists a constant $L_D$ such that for any $x, y \in \Omega$ and for any $g_u \in \partial D(\cdot, x)(y)$, there exists $g_r \in \partial r(x)$ such that $\|g_u - g_r\| \leq L_D\|x - y\|$.
	\end{assumption}
	Now we present the convergence rate results for our method.
	\begin{theorem}\label{the3.19}
		\textbf
		Let $\{x_k\}$ be the sequence generated by Algorithm \ref{PSGA}. Suppose Assumptions \ref{ass3.4}-\ref{ass3.7} hold. Then
		\begin{align}\notag
			\min_{k=1,2,...K}  \mathbb{E}\textup{dist}(0, \partial F(x_k))\leq \sqrt{\frac{G_3}{K}} = \mathcal{O}\left(\sqrt{\frac{1}{K}}\right),
		\end{align}\textup{where} \[
		\textup{dist}\left(0, \partial F(x^k)\right) := \inf_{v \in \partial F(x^k)} \|v\|.
		\]
	\end{theorem}
	
	\noindent  Proof. From (\ref{3.28}) we have
	\begin{align*}
		g_u := -\frac{y_k-x_k}{\eta_k} - \widetilde{\nabla}f(x_k) \in \partial D(\cdot, x_k)(y_k),
	\end{align*}
	which  by Assumption \ref{ass3.11} implies that there exists $ g_r \in \partial r(x_k)$ such that $\| g_r - g_u\| \leq L_D\|y_k - x_k\|.$ Therefore,
	\begin{align}
		\textup{dist}(0, \partial F(x_k))
		&\leq \| \nabla f(x_k) + g_r \| \notag \\
		& = \| \nabla f(x_k)  -\frac{y_k-x_k}{\eta_k} - \widetilde{\nabla}f(x_k) + g_r - g_u \| \notag\\
		& \leq \| \nabla f(x_k) - \widetilde{\nabla}f(x_k)\| + \|\frac{y_k-x_k}{\eta_k}\| + \|g_r - g_u\| \notag\\
		& \leq \| \nabla f(x_k) - \widetilde{\nabla}f(x_k) \| +(\frac{1}{\eta_k}+L_D)\|y_k-x_k\|.\label{3.52}
	\end{align}

	Taking the total expectation in (\ref{3.52}), we get
	\begin{align}
		\mathbb{E} \textup{dist}^2(0,\partial F(x_k))
		&\leq 2 \mathbb{E}\| \nabla f(x_k) - \widetilde{\nabla}f(x_k) \|^2+2(\frac{1}{\eta_k}+L_D)^2\mathbb{E}\|y_k-x_k\|^2 \notag \\
		& =  2\mathbb{E}\| \nabla f(x_k) - \widetilde{\nabla}f(x_k) \|^2 + 2(\frac{1}{\eta_k}+L_D)^2 \frac{(k+1)^2}{k^2}\mathbb{E} \|x_{k+1}-x_k\|^2\notag\\
		& \leq 2\mathbb{E}\| \nabla f(x_k) - \widetilde{\nabla}f(x_k) \|^2 +8(\frac{1}{C_0}+L_D)^2\mathbb{E}\|x_{k+1}-x_{k}\|^2, \label{3.56}
	\end{align}
	where the first inequality follows from the fact that $\|a+b\|^2\leq2\|a\|^2+2\|b\|^2,\forall a,b\in \mathbb{R}^n$ and the second one is due to (\ref{in3.20}).
	
	We sum the inequality (\ref{3.56}) for $k=1,2,\cdots,K$ to obtain
	\begin{align}
		\sum_{k=1}^{K} \mathbb{E}\textup{dist}^2(0,\partial F(x_k))
		&\leq 2\sum_{k=1}^{K}\mathbb{E}\| \nabla f(x_k) - \widetilde{\nabla}f(x_k) \|^2+8(\frac{1}{C_0}+L_D)^2\sum_{k=1}^{K}\mathbb{E}\|x_{k+1}-x_{k}\|^2 \notag\\
		&\leq 2G_2+ 8(\frac{1}{C_0}+L_D)^2G_1,
	\end{align}
	where the second inequality follows from (\ref{3.14}) and (\ref{in3.30}).
	
	Setting $G_3 = 2G_2+ 8(\frac{1}{C_0}+L_D)^2G_1$, we get
	\begin{equation}
		\min_{k=1,2,...K}  \mathbb{E}\textup{dist}^2(0, \partial F(x_k)) \leq \frac{G_3}{K},
	\end{equation}
	which completes the proof.

	\section{Numerical Experiments}\label{Section4}
	
	In this section, we analyze the efficiency of our PSGA algorithm and compare it with other algorithms employing variance reduction techniques. The comparison focuses on two aspects: convergence rates and gradient estimation errors. All experiments were run on a computer with an AMD Ryzen 7 5800H 3.20 GHz CPU and 16GB of memory.

	We evaluate the algorithms on two standard problems: Logistic regression with $\ell_1$-regularization and Lasso regression. We compare our Algorithm \ref{PSGA} (PSGA) with S-PStorm\cite{dai2023variance}, SAGA\cite{defazio2014saga}, RDA\cite{xiao2009dual}, Prox-SVRG\cite{xiao2014proximal}, and PStorm\cite{xu2023momentum} algorithms.
	
	The parameters of each algorithm are set as follows:
	
	$(a)$
	For ProxSVRG, SAGA, and S-PStorm algorithms, we use a constant step size strategy by setting $\alpha_{k}\equiv0.1/L$. For RDA algorithm we set step size as $\eta_k=\sqrt{k}/\gamma$, where $\gamma=10^{-2}$ is suggested in \cite{dai2023variance}. For Pstorm algorithm we take $\eta_{k}=\dfrac{4^{1/3}/8L}{(k+4)^{1/3}}$ as in \cite{xu2023momentum}.
	
	$(b)$
	For PStorm algorithm we take $\beta_{k} = \dfrac{1+24\eta_{k}^2L^2-\frac{\eta_{k+1}}{\eta_{k}}}{1+4\eta_{k}^2L^2}$. For S-PStorm algorithm  we take $\beta_k = \dfrac{1}{k+1}$. For  our algorithm(PSGA) we take $\theta_{k}=\dfrac{1}{k+1}$.
	
	In our numerical experiments, we imposed a stopping rule: a test is terminated when either the maximum number of 1000 iterations is reached or the 12-hour runtime limit is reached.
	
	\begin{table}[ht]
		\centering
		\caption{Datasets used in experiments}
		\label{table1}
		\begin{tabularx}{0.9\textwidth}{|l|>{\centering\arraybackslash}X|>{\centering\arraybackslash}X|}
			\hline
			\textbf{dataset} & \textbf{Data Points Number N} & \textbf{Feature Number n} \\
			\hline
			a9a & 32,561 & 123 \\
			covtype & 581,012 & 54 \\
			phishing & 11,055 & 68 \\
			rcv1 & 20,242 & 47,236 \\
			real-sim & 72,309 & 20,958 \\
			news20 & 19,996 & 1,355,191 \\
			w8a & 49,749 & 300 \\
			\hline
		\end{tabularx}
	\end{table}
	Datasets for Logistic regression with $\ell_1$-regularization and Lasso regression problems are obtained from the LIBSVM \cite{chang2011libsvm}. Details of the datasets and the parameters are given in Table \ref{table1}.

	\subsection{Logistic Regression Problem}
	We consider solving problem (\ref{1.1}) given by the regularized binary Logistic loss with $L$-smooth convex function and non-smooth convex group-$\ell_1$ regularizer:
	\[
	\min_{x \in \mathbb{R}^n} \;\;
	\frac{1}{N} \sum_{j=1}^N \log\bigl(1 + e^{-y_j\, x^T d_j}\bigr)
	\;+\; 10^{-5} \,\|x\|_1^2
	\;,
	\]
	where \( N \) is the number of data points, \( d_j \in \mathbb{R}^n \) is the \( j \)-th data point, and \( y_j \in \{-1, 1\} \) is the class label for the \( j \)-th data point. In the following figures, $f^*$ represents the lowest objective function value obtained among all tested algorithms.
	
	\begin{figure}[htbp]
		\centering

		\begin{minipage}[t]{0.48\textwidth}
			\centering
			\vspace{0pt}
			{\footnotesize\textbf{a9a}}
			\vspace{1pt}
			\includegraphics[width=\linewidth, height=4.5cm, keepaspectratio]{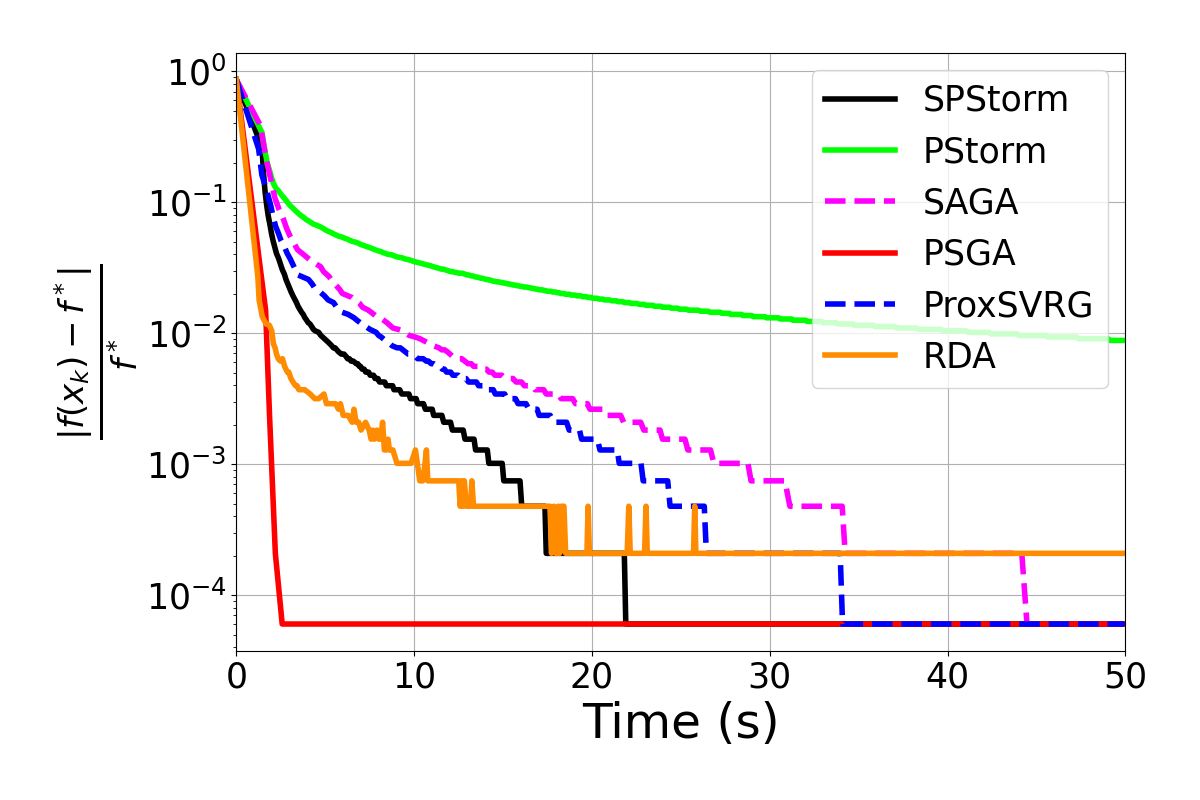}
			\label{fig:fx1}
		\end{minipage}
		\hfill
		\begin{minipage}[t]{0.48\textwidth}
			\centering
			\vspace{0pt}
			{\footnotesize\textbf{covtype}}
			\vspace{1pt}
			\includegraphics[width=\linewidth, height=4.5cm, keepaspectratio]{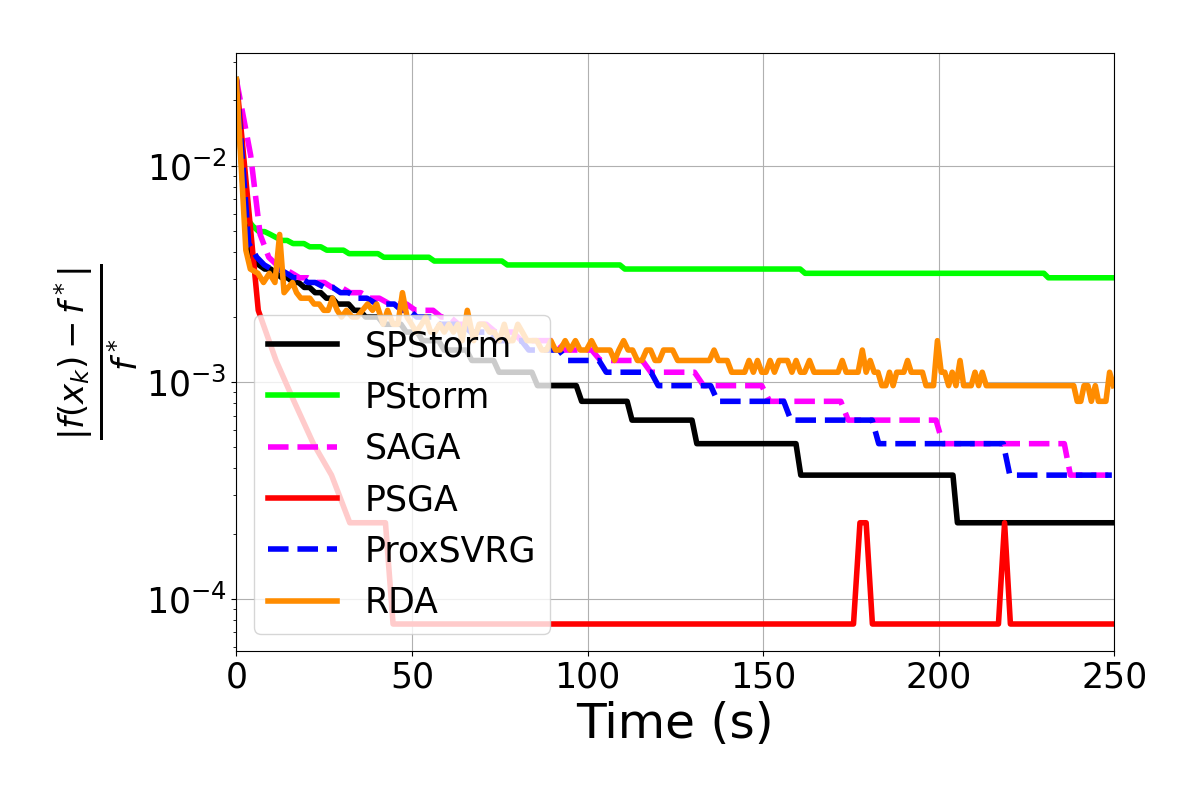}
			\label{fig:fx2}
		\end{minipage}
		
		\vspace{4pt}

		\begin{minipage}[t]{0.48\textwidth}
			\centering
			\vspace{0pt}
			{\footnotesize\textbf{phishing}}
			\vspace{1pt}
			\includegraphics[width=\linewidth, height=4.5cm, keepaspectratio]{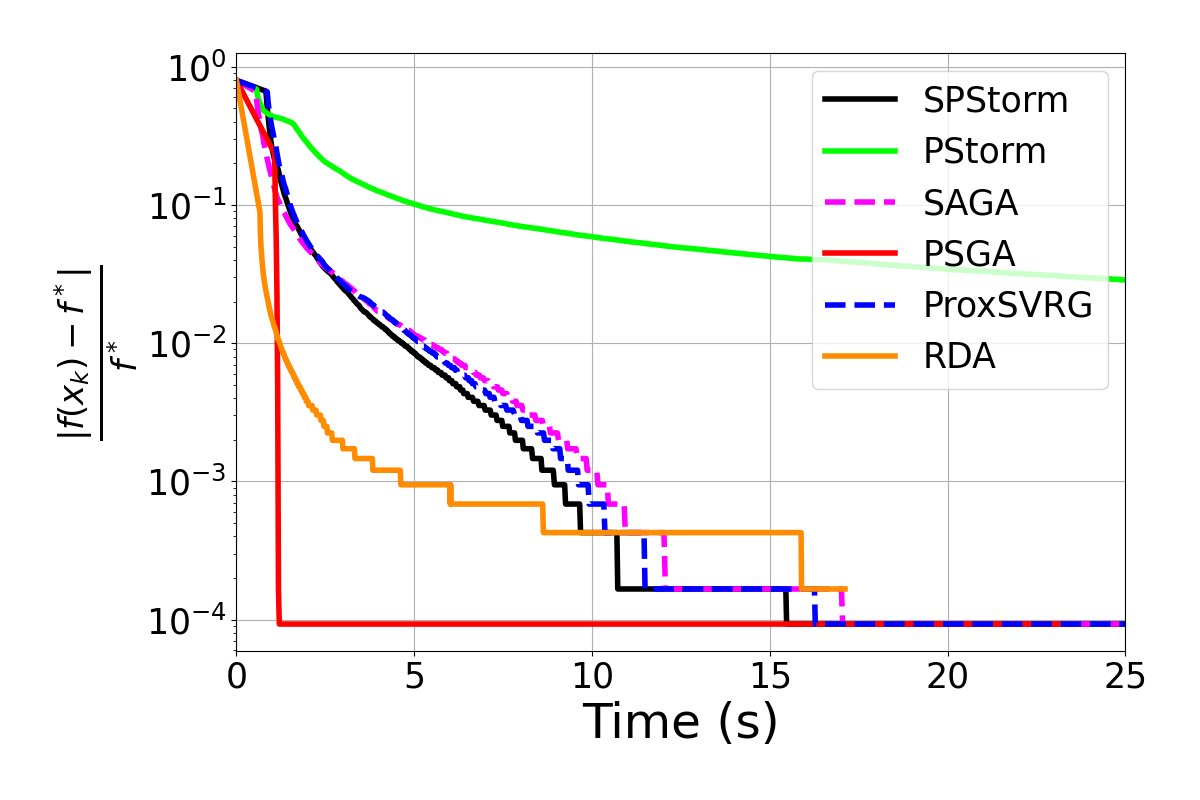}
			\label{fig:fx3}
		\end{minipage}
		\hfill
		\begin{minipage}[t]{0.48\textwidth}
			\centering
			\vspace{0pt}
			{\footnotesize\textbf{rcv1}}
			\vspace{1pt}
			\includegraphics[width=\linewidth, height=4.5cm, keepaspectratio]{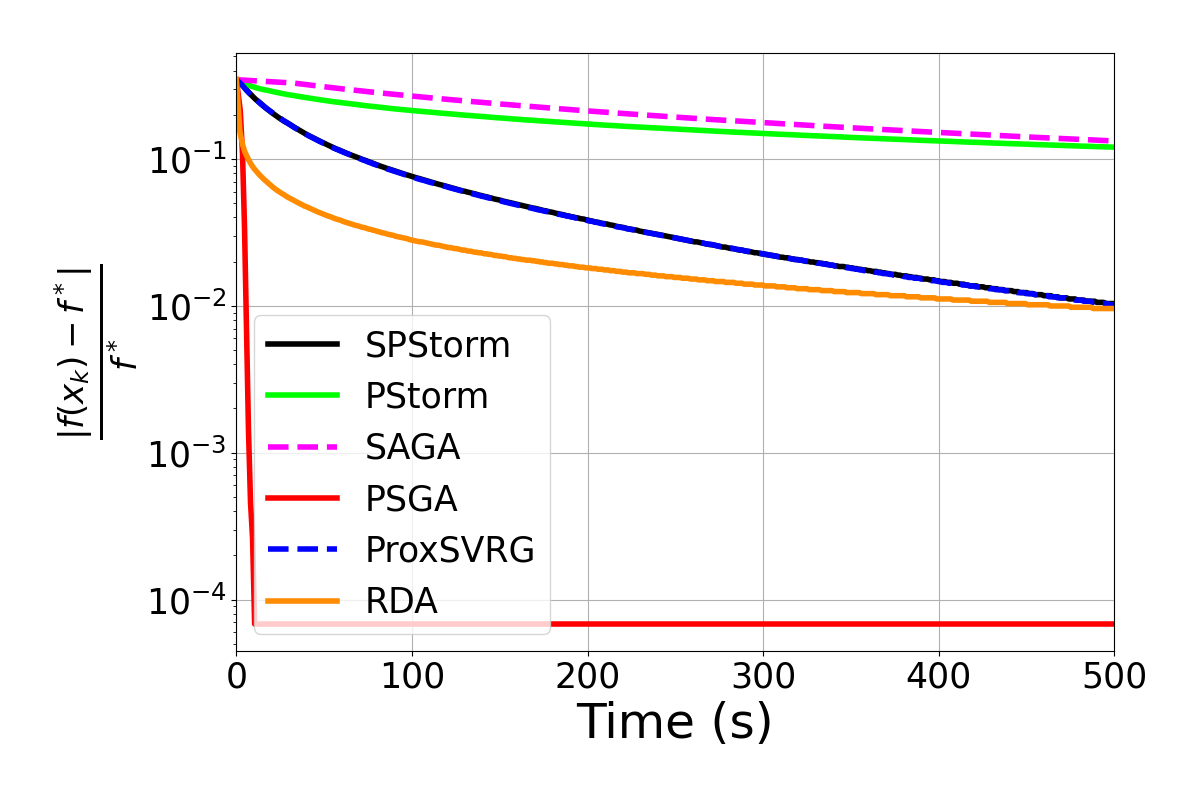}
			\label{fig:fx4}
		\end{minipage}
		
		\vspace{4pt}

		\begin{minipage}[t]{0.48\textwidth}
			\centering
			\vspace{0pt}
			{\footnotesize\textbf{real-sim}}
			\vspace{1pt}
			\includegraphics[width=\linewidth, height=4.5cm, keepaspectratio]{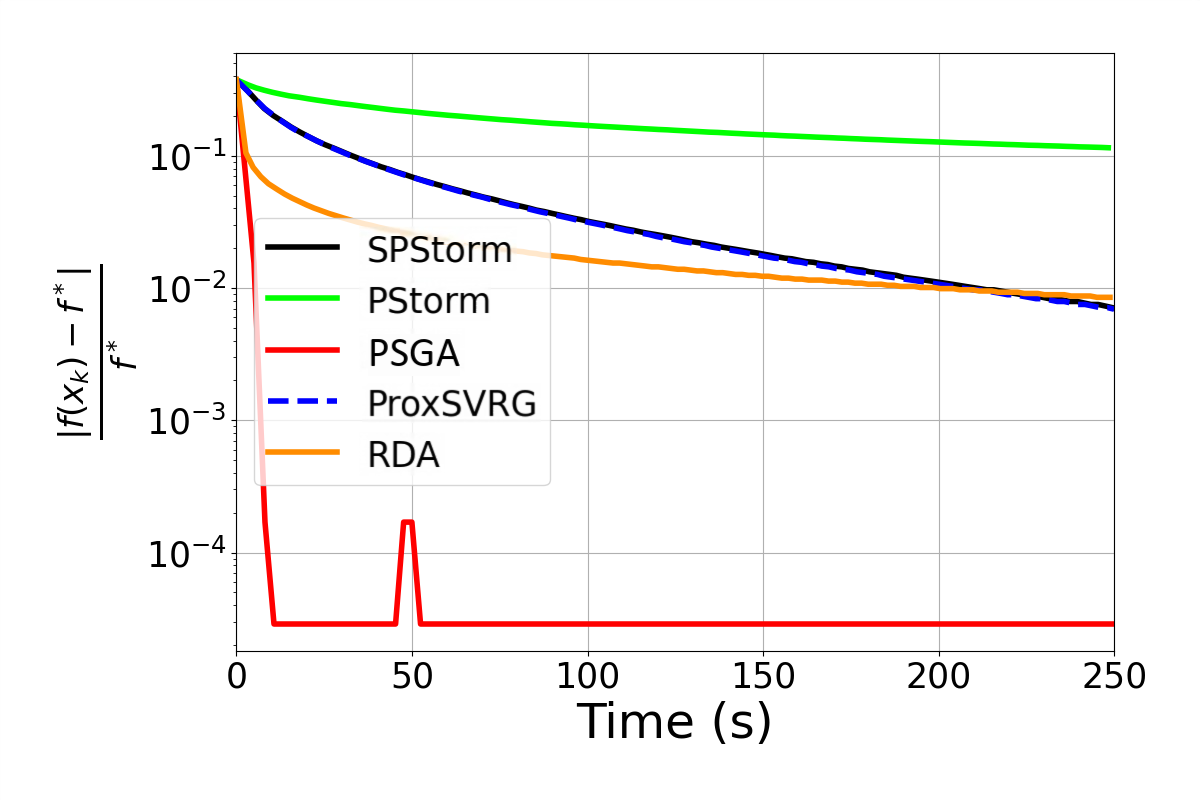}
			\label{fig:fx5}
		\end{minipage}
		\hfill
		\begin{minipage}[t]{0.48\textwidth}
			\centering
			\vspace{0pt}
			{\footnotesize\textbf{w8a}}
			\vspace{1pt}
			\includegraphics[width=\linewidth, height=4.5cm, keepaspectratio]{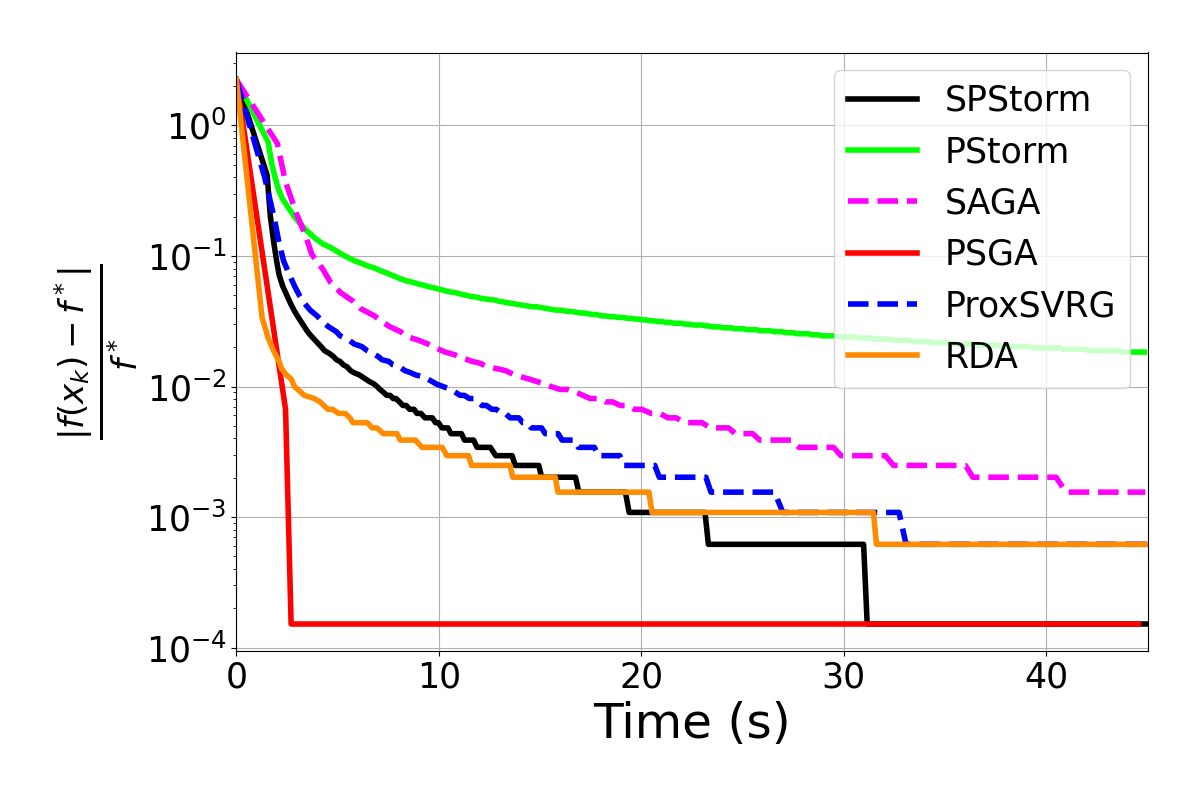}
			\label{fig:fx6}
		\end{minipage}
		
		\vspace{4pt}

		\begin{minipage}[t]{0.48\textwidth}
			\centering
			\vspace{0pt}
			{\footnotesize\textbf{news20}}
			\vspace{1pt}
			\includegraphics[width=\linewidth, height=4.5cm, keepaspectratio]{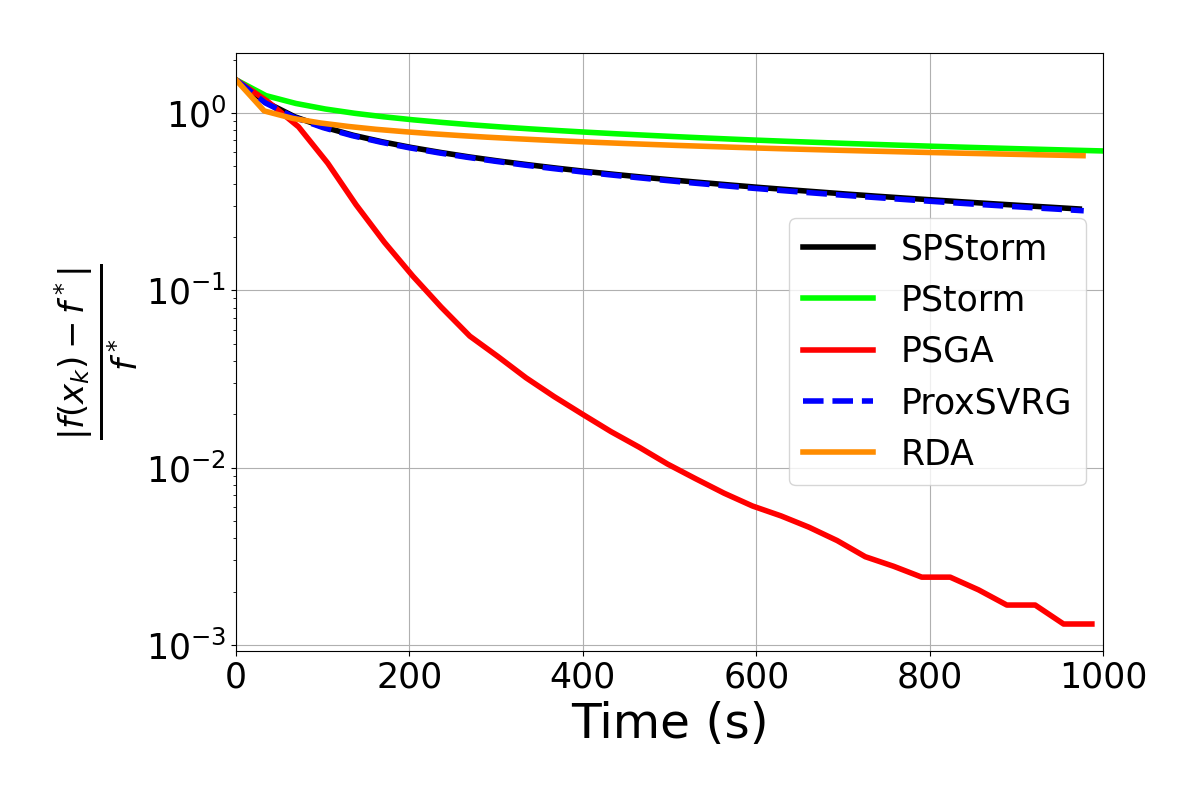}
			\label{fig:fx7}
		\end{minipage}
		
		\vspace{4pt}
		\caption{Evolution of $\frac{|f(x_k) - f^*|}{f^*}$ with respect to runtime on a9a, covtype, phishing, rcv1, real-sim, news20 and w8a. }
		\label{fig:fx_results}
	\end{figure}
	
	\begin{figure}[htbp]
		\centering

		\begin{minipage}[t]{0.48\textwidth}
			\centering
			\vspace{0pt}
			{\footnotesize\textbf{a9a}}
			\vspace{1pt}
			\includegraphics[width=\linewidth, height=4.5cm, keepaspectratio]{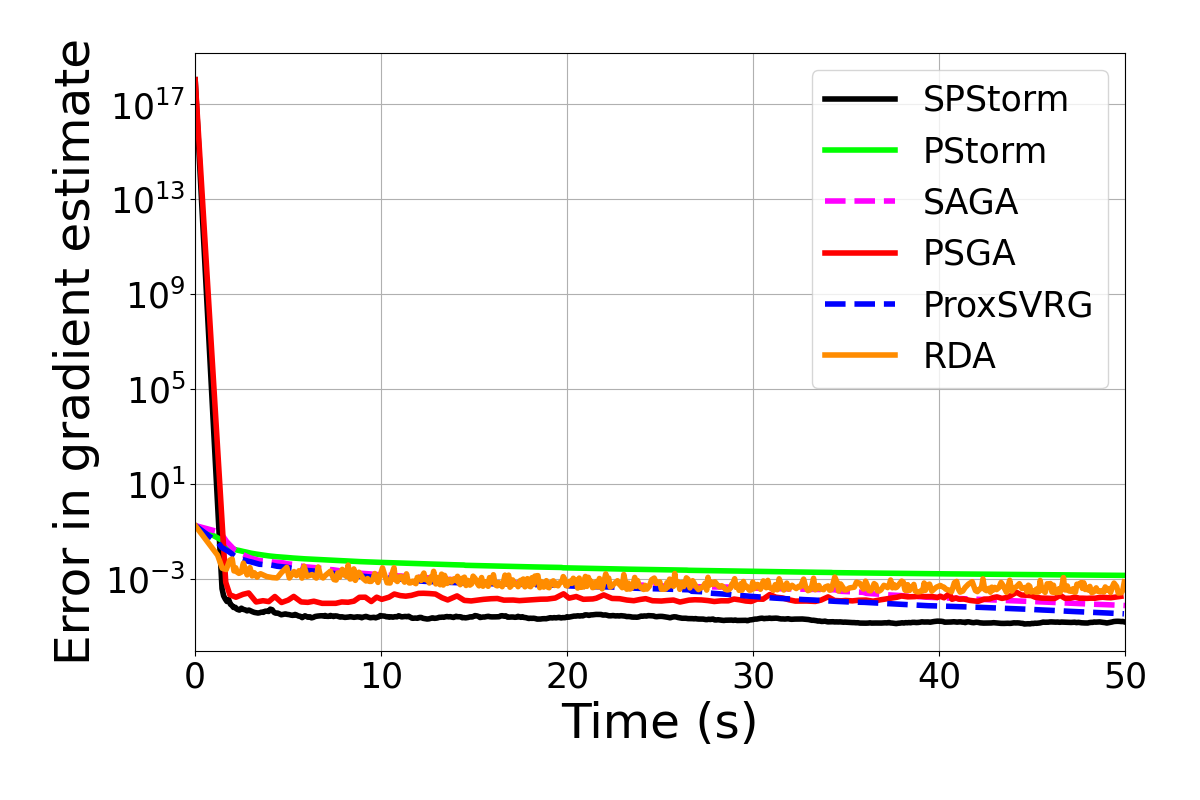}
			\label{fig:er1}
		\end{minipage}
		\hfill
		\begin{minipage}[t]{0.48\textwidth}
			\centering
			\vspace{0pt}
			{\footnotesize\textbf{covtype}}
			\vspace{1pt}
			\includegraphics[width=\linewidth, height=4.5cm, keepaspectratio]{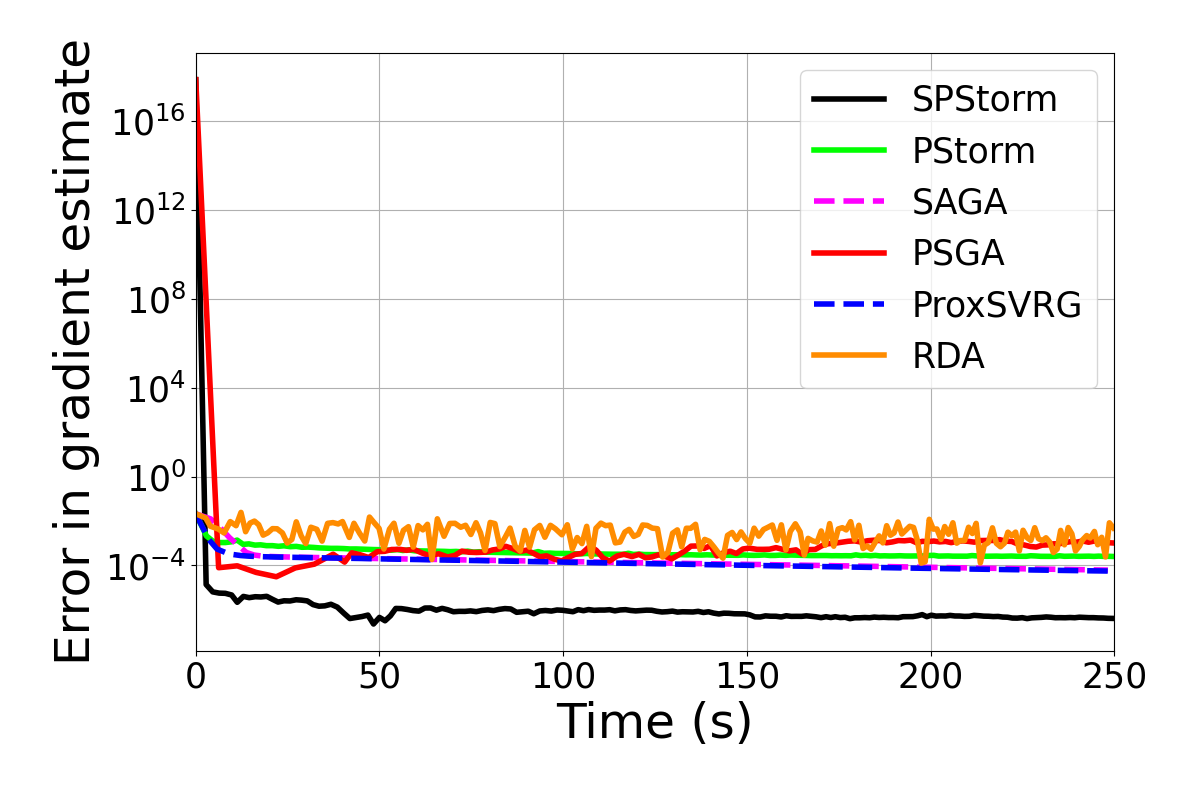}
			\label{fig:er2}
		\end{minipage}
		
		\vspace{4pt}

		\begin{minipage}[t]{0.48\textwidth}
			\centering
			\vspace{0pt}
			{\footnotesize\textbf{phishing}}
			\vspace{1pt}
			\includegraphics[width=\linewidth, height=4.5cm, keepaspectratio]{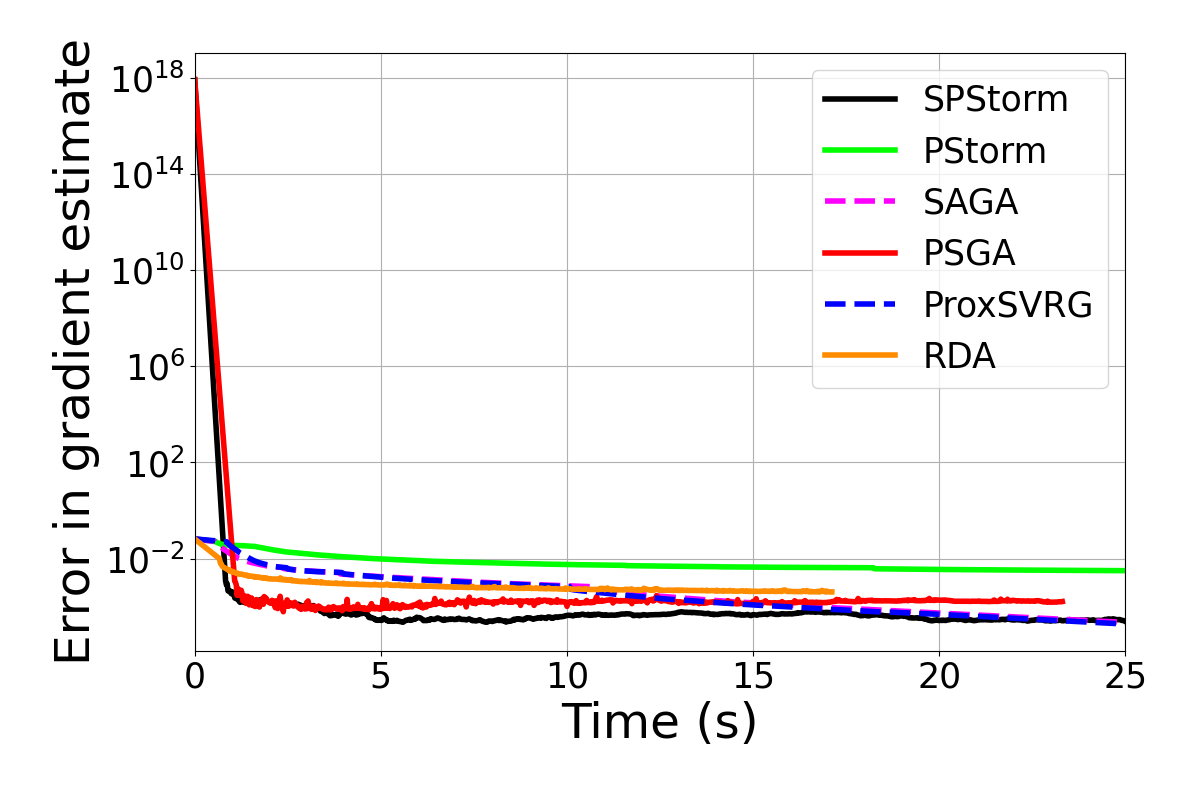}
			\label{fig:er3}
		\end{minipage}
		\hfill
		\begin{minipage}[t]{0.48\textwidth}
			\centering
			\vspace{0pt}
			{\footnotesize\textbf{rcv1}}
			\vspace{1pt}
			\includegraphics[width=\linewidth, height=4.5cm, keepaspectratio]{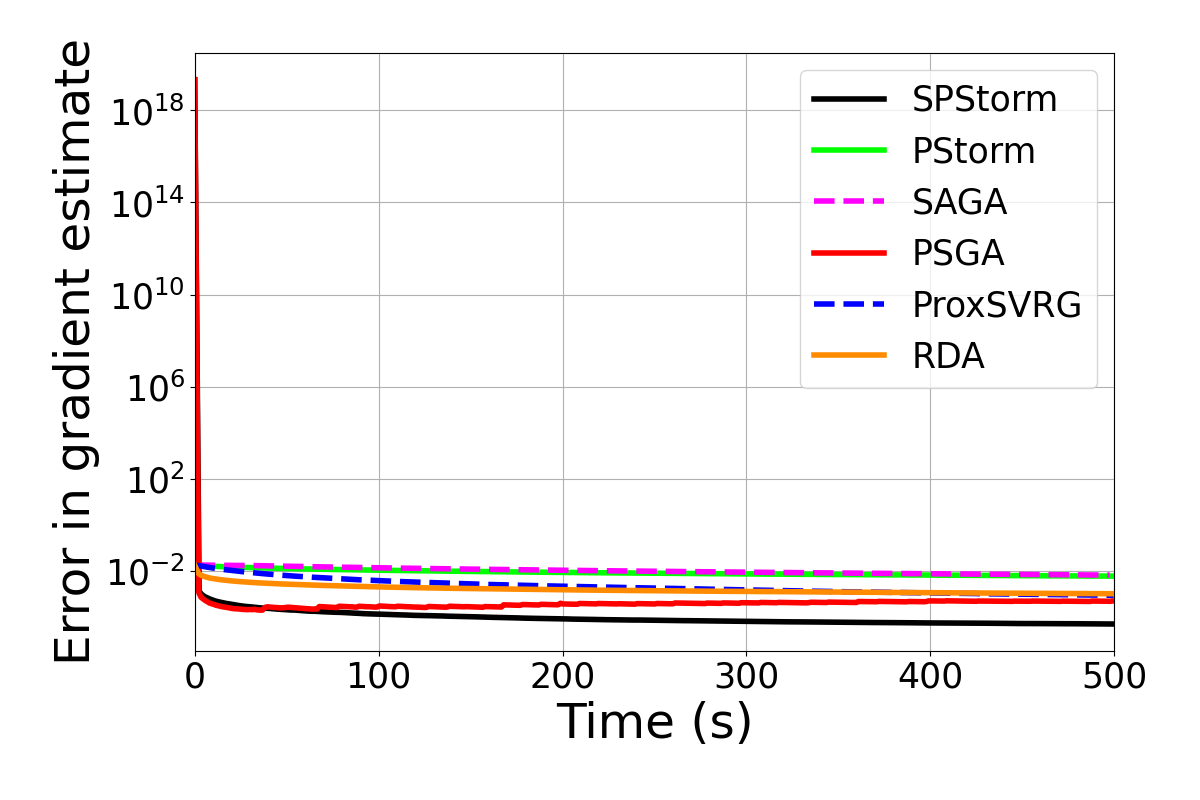}
			\label{fig:er4}
		\end{minipage}
		
		\vspace{4pt}

		\begin{minipage}[t]{0.48\textwidth}
			\centering
			\vspace{0pt}
			{\footnotesize\textbf{real-sim}}
			\vspace{1pt}
			\includegraphics[width=\linewidth, height=4.5cm, keepaspectratio]{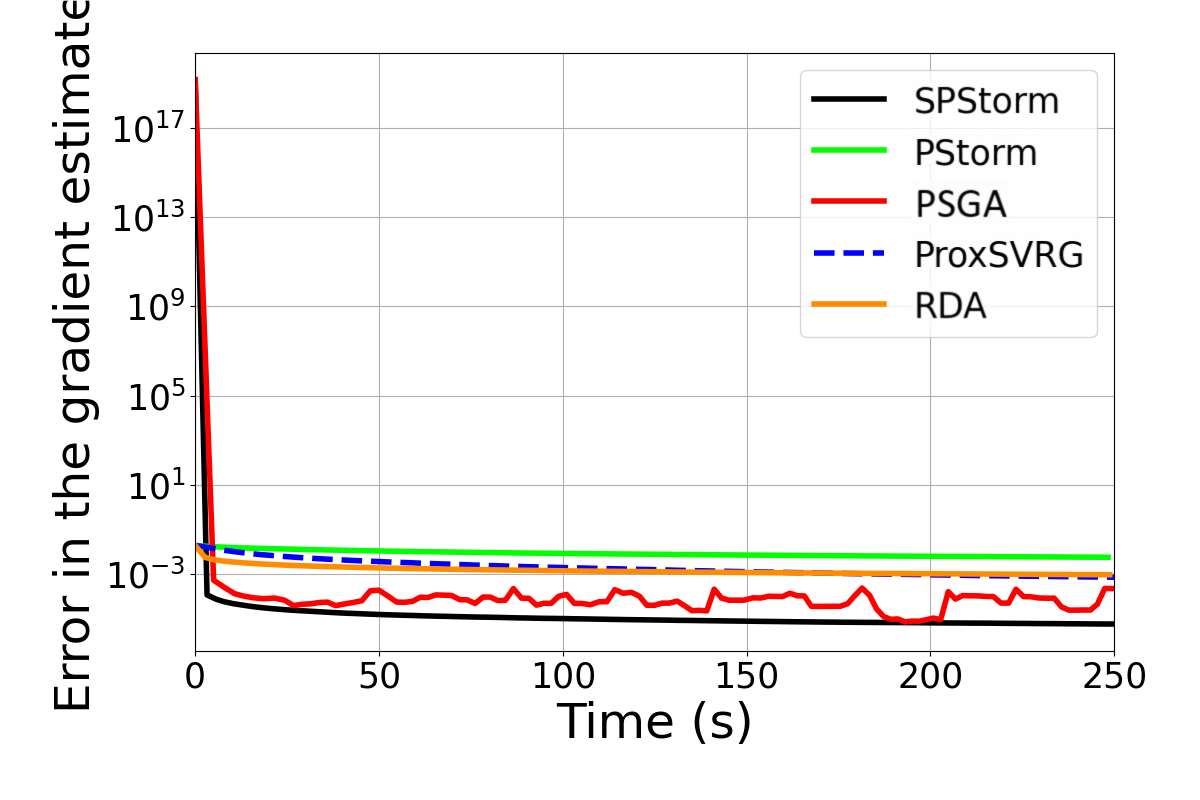}
			\label{fig:er5}
		\end{minipage}
		\hfill
		\begin{minipage}[t]{0.48\textwidth}
			\centering
			\vspace{0pt}
			{\footnotesize\textbf{w8a}}
			\vspace{1pt}
			\includegraphics[width=\linewidth, height=4.5cm, keepaspectratio]{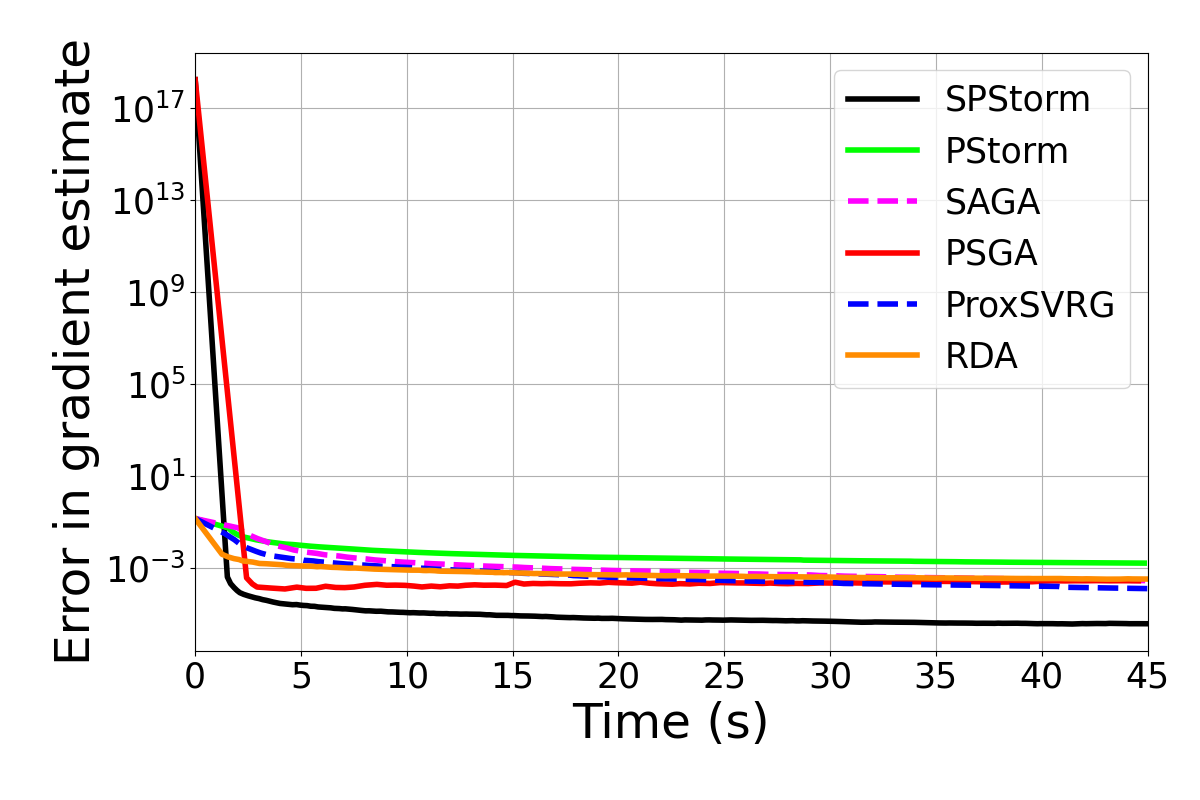}
			\label{fig:er6}
		\end{minipage}
		
		\vspace{4pt}

		\begin{minipage}[t]{0.48\textwidth}
			\centering
			\vspace{0pt}
			{\footnotesize\textbf{news20}}
			\vspace{1pt}
			\includegraphics[width=\linewidth, height=4.5cm, keepaspectratio]{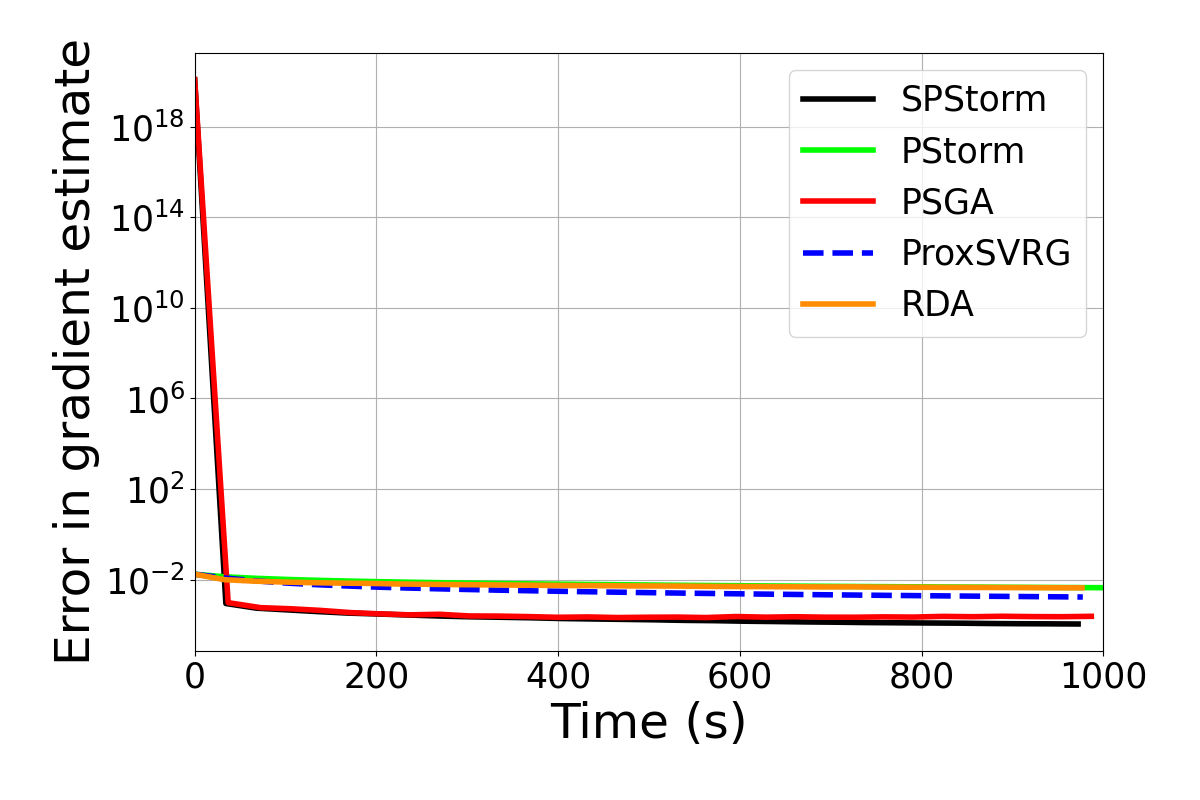}
			\label{fig:er7}
		\end{minipage}
		
		\vspace{4pt}
		\caption{ Evolution of gradient estimation error with respect to runtime on a9a, covtype, phishing, rcv1, real-sim, news20 and w8a.}
		\label{fig:er_results}
	\end{figure}
	
	In Figure \ref{fig:fx_results}, we observe that our algorithm(PSGA) achieves faster convergence across all datasets.

	From Figure \ref{fig:er_results}, we can see that our algorithm(PSGA) has smaller gradient estimation error than other five methods on the datasets phishing, rcv1 and news20, and hence our method has higher accuracy. For datasets a9a and real-sim, we find that the gradient estimation errors of S-PStorm method are almost the same with ours, but our method needs fewer CPU time.
	
	Table \ref{table2}  presents the minimum values $f(best)$ achieved by each method, along with the computation time and the number of iterations required to reach the $f(best)$. The symbol ``$-$" indicates that the algorithm cannot be tested on the data set.
	
	\begin{table}[htbp]
		\centering
		\caption{Convergence Performance on Different Datasets}
		\label{table2}

		\begin{subtable}{0.48\textwidth}
			\centering
			\caption{a9a dataset}
			\begin{tabular}{l S[table-format=1.6] S[table-format=4.0] S[table-format=5.2]}
				\toprule
				Algorithm & {$f(best)$} & {Iter.} & {Time (s)} \\
				\midrule
				{PSGA}      & {0.3723}      &    {4}    &  {2.57}       \\
				PStorm    & 0.3742      &  968    & 94.38      \\
				ProxSVRG  & 0.3723      &  217    & 34.09      \\
				RDA       & 0.3723      &  721    & 71.98      \\
				SAGA      & 0.3723      &  218    & 44.47      \\
				SPStorm   & 0.3723      &  217    & 21.90      \\
				\bottomrule
			\end{tabular}
		\end{subtable}
		\hfill
		\begin{subtable}{0.48\textwidth}
			\centering
			\caption{covtype dataset}
			\begin{tabular}{l S[table-format=1.6] S[table-format=4.0] S[table-format=5.2]}
				\toprule
				Algorithm & {$f(best)$} & {Iter.} & {Time (s)} \\
				\midrule
				{PSGA}      & {0.6762}      &   {79}    &  {175.73}     \\
				PStorm     & 0.6776      &  950    & 1287.68    \\
				ProxSVRG   & 0.6762      &  662    & 1057.77    \\
				RDA       & 0.6765      &  661    & 678.50     \\
				SAGA      & 0.6762      &  663    & 1083.31    \\
				SPStorm    & 0.6762      &  662    & 883.59     \\
				\bottomrule
			\end{tabular}
		\end{subtable}

		\begin{subtable}{0.48\textwidth}
			\centering
			\caption{phishing dataset}
			\begin{tabular}{l S[table-format=1.6] S[table-format=4.0] S[table-format=5.2]}
				\toprule
				Algorithm & {$f(best)$} & {Iter.} & {Time (s)} \\
				\midrule
				{PSGA}      & {0.3857}      &   {7}    &  {1.20}      \\
				PStorm     & 0.3957      &  999    & 27.85      \\
				ProxSVRG   & 0.3857      &  556    & 16.28      \\
				RDA       & 0.3858      &  927    & 15.89      \\
				SAGA      & 0.3857      &  557    & 17.05      \\
				SPStorm    & 0.3857      &  553    & 15.47      \\
				\bottomrule
			\end{tabular}
		\end{subtable}
		\hfill
		\begin{subtable}{0.48\textwidth}
			\centering
			\caption{rcv1 dataset}
			\begin{tabular}{l S[table-format=1.6] S[table-format=4.0] S[table-format=5.2]}
				\toprule
				Algorithm & {$f(best)$} & {Iter.} & {Time (s)} \\
				\midrule
				{PSGA}      & {0.5148}      &   {12}    & {20.80}      \\
				PStorm     & 0.5549      &  999    & 1210.62    \\
				ProxSVRG   & 0.5155      &  963    & 1185.14    \\
				RDA       & 0.5173      &  967    & 1107.71    \\
				SAGA      & 0.5515      &  963    & 15635.52   \\
				SPStorm    & 0.5155      &  963    & 1179.82    \\
				\bottomrule
			\end{tabular}
		\end{subtable}

		\begin{subtable}{0.48\textwidth}
			\centering
			\caption{real-sim dataset}
			\begin{tabular}{l S[table-format=1.6] S[table-format=4.0] S[table-format=5.2]}
				\toprule
				Algorithm & {$f(best)$} & {Iter.} & {Time (s)} \\
				\midrule
				{PSGA}      & {0.5035}      &    {4}    & {10.71}      \\
				PStorm     & 0.5190      & 1001    & 1902.18    \\
				ProxSVRG   & 0.5035      &  510    & 1017.24    \\
				RDA       & 0.5043      &  981    & 1733.21    \\
				SAGA      & \multicolumn{3}{c}{--}           \\
				SPStorm    & 0.5035      &  510    & 1027.93    \\
				\bottomrule
			\end{tabular}
		\end{subtable}
		\hfill
		\begin{subtable}{0.48\textwidth}
			\centering
			\caption{news20 dataset}
			\begin{tabular}{l S[table-format=1.6] S[table-format=4.0] S[table-format=5.2]}
				\toprule
				Algorithm & {$f(best)$} & {Iter.} & {Time (s)} \\
				\midrule
				{PSGA}      & {0.2727}      &  {64}    &  {2128.66}    \\
				PStorm     & 0.3152      & 1000    & 33511.46   \\
				ProxSVRG   & 0.2739      &  982    & 32190.78   \\
				RDA       & 0.3354      & 1000    & 31708.78   \\
				SAGA      & \multicolumn{3}{c}{--}           \\
				SPStorm    & 0.2729      &  982    & 35990.36   \\
				\bottomrule
			\end{tabular}
		\end{subtable}

		\begin{subtable}{\textwidth}
			\centering
			\caption{w8a dataset}
			\begin{tabular}{l S[table-format=1.6] S[table-format=4.0] S[table-format=5.2]}
				\toprule
				Algorithm & {$f(best)$} & {Iter.} & {Time (s)} \\
				\midrule
				{PSGA}      & {0.4265}      &    {10}    &  {2.10}      \\
				PStorm     & 0.4265      &  629    & 78.07      \\
				ProxSVRG   & 0.4265      &   39    &  5.72      \\
				RDA       & 0.4265      &   80    &  7.02      \\
				SAGA      & 0.4265      &   40    &  9.87      \\
				SPStorm    & 0.4265      &   39    &  5.38      \\
				\bottomrule
			\end{tabular}
		\end{subtable}
	\end{table}
	
	From Table \ref{table2}, it can be observed that our algorithm(PSGA) obtains objective function values $f(best)$ that is no worse than those of other algorithms across all tested datasets. At the same time, our algorithm requires fewer iterations and less CPU time than other algorithms. Additionally, we note that SAGA terminated immediately on the datasets news20 and real-sim because the storage of the gradient look-up table exceeded the memory limit.
	
	\newpage
	\subsection{Lasso Regression Problem}
	
	In this section we consider solving problem (\ref{1.1}) given by the Lasso loss with $L$-smooth convex function and non-smooth convex group-$\ell_1$ regularizer:
	
	\[\min_{x \in \mathbb{R}^n} \;\;
	\frac{1}{2N} \sum_{i=1}^n (y_i - \mathbf{A}^T \mathbf{x})^2
	\;+\; 10^{-5} \,\|x\|_1^2
	\;.\]
	where $A$ is characteristic matrix and $y_i$ is the true value of the sample. The following figures and table are our experiment results:
	\begin{figure}[htbp]
		\centering

		\begin{minipage}[t]{0.48\textwidth}
			\centering
			\textbf{news20}\\
			\vspace{0.2em}
			\includegraphics[width=\linewidth, height=5cm, keepaspectratio]{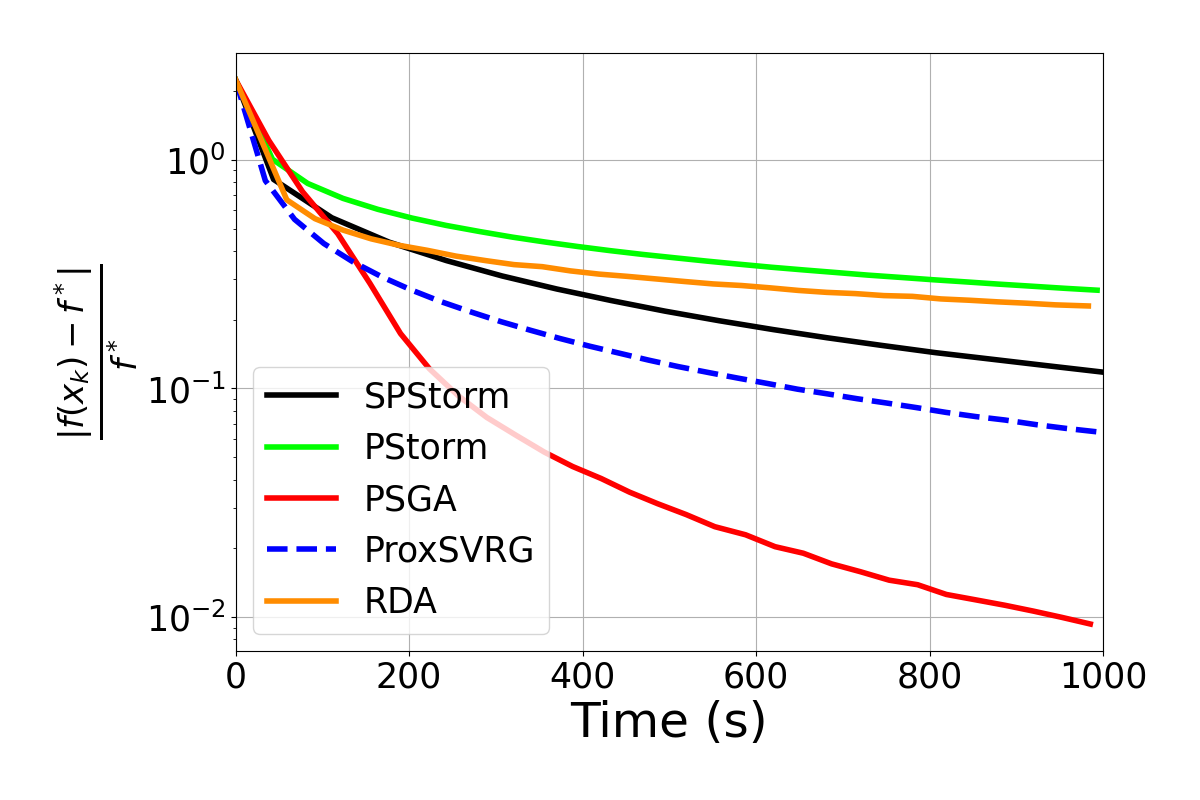}
			\label{fig:fx_news20}
		\end{minipage}
		\hfill
		\begin{minipage}[t]{0.48\textwidth}
			\centering
			\textbf{rcv1}\\
			\vspace{0.2em}
			\includegraphics[width=\linewidth, height=5cm, keepaspectratio]{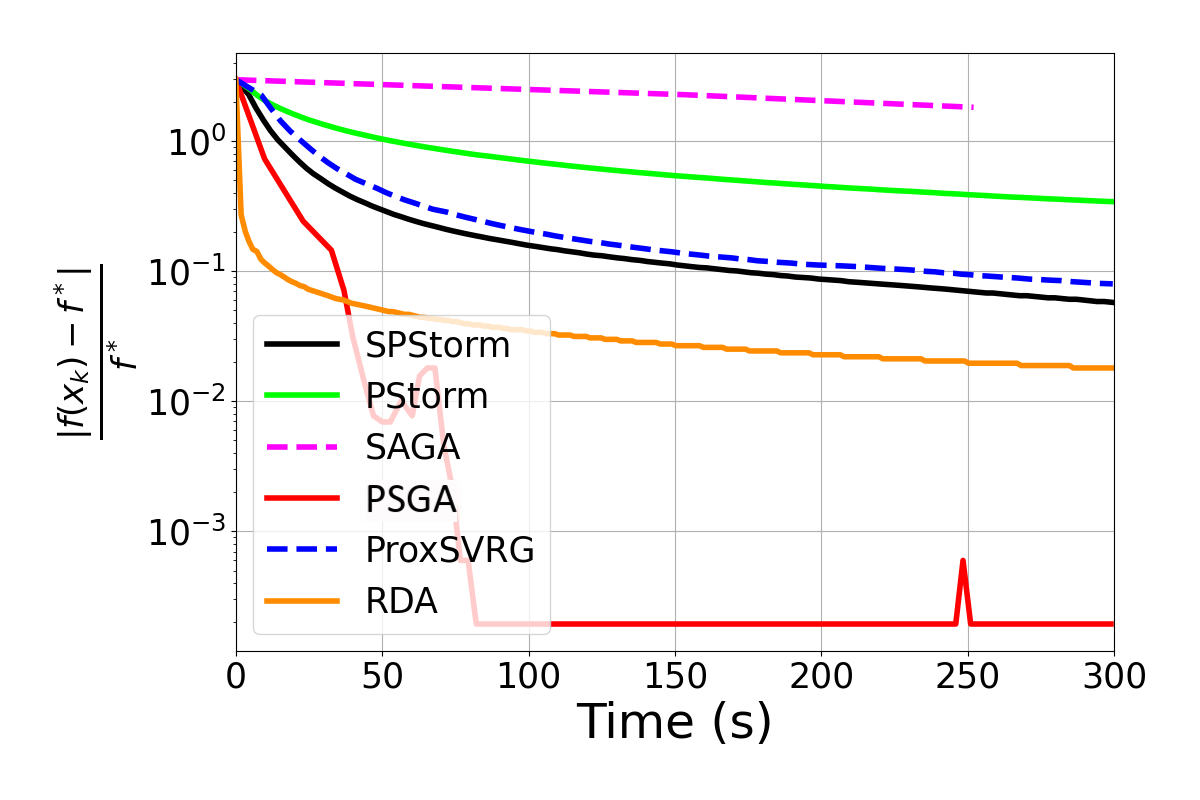}
			\label{fig:fx_rcv1}
		\end{minipage}
		
		\vspace{0.5em}
		\caption{\text{Evolution of $\frac{|f(x_k) - f^*|}{f^*}$ with respect to runtime on rcv1 and news20. }}
		\label{fig:lafx_results}
	\end{figure}
	
	\begin{figure}[htbp]
		\centering

		\begin{minipage}[t]{0.48\textwidth}
			\centering
			\textbf{news20}\\
			\vspace{0.2em}
			\includegraphics[width=\linewidth, height=5cm, keepaspectratio]{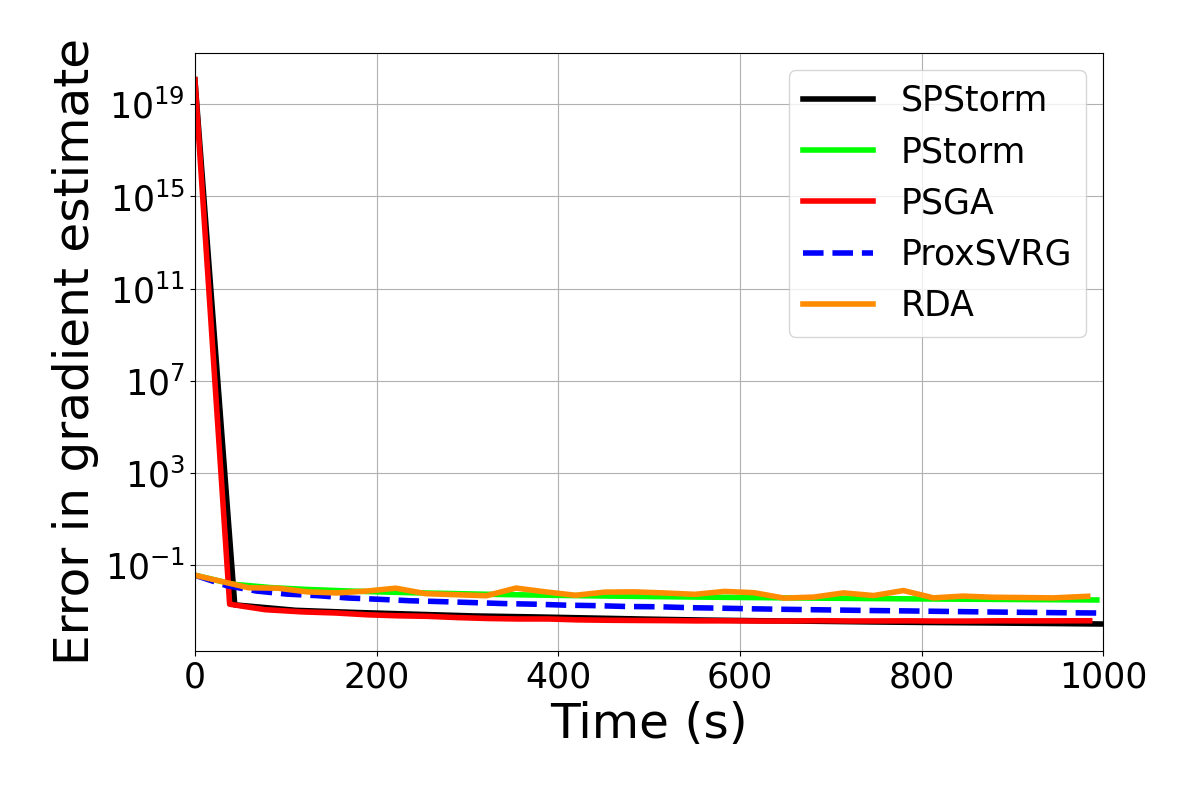}
			\label{fig:er_news20}
		\end{minipage}
		\hfill
		\begin{minipage}[t]{0.48\textwidth}
			\centering
			\textbf{rcv1}\\
			\vspace{0.2em}
			\includegraphics[width=\linewidth, height=5cm, keepaspectratio]{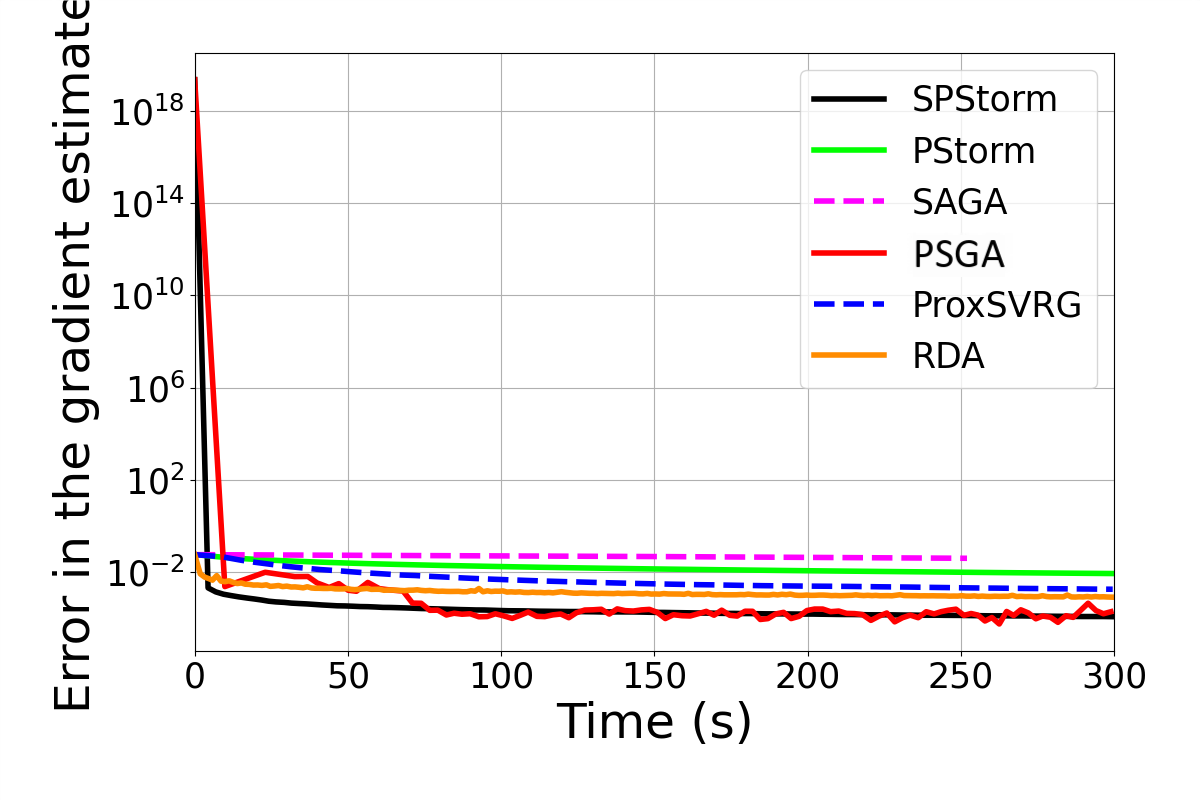}
			\label{fig:er_rcv1}
		\end{minipage}
		
		\vspace{0.5em}
		\caption{\text{ Evolution of gradient estimation error with respect to runtime on rcv1 and news20 }}
		\label{fig:laer_results}
	\end{figure}
	
	From Figures \ref{fig:lafx_results} and \ref{fig:laer_results}, we observe that our algorithm(PSGA) achieves faster convergence and achieves more precise gradient estimates on the rcv1 and news20 datasets.
	
	\begin{table}[htbp]
		\centering
		\caption{Convergence Performance on Different Datasets}
		\label{tab:convergence_all}

		\begin{subtable}{0.48\textwidth}
			\centering
			\caption{news20 dataset}
			\begin{tabular}{l S[table-format=1.6] S[table-format=4.0] S[table-format=5.1]}
				\toprule
				Algorithm & {$f$(best)} & {Iter.} & {Time (s)} \\
				\midrule
				{PSGA}      & {0.1547}      &  {99}    & {3116.1}    \\
				PStorm    & 0.1580      &  971    & 36416.9    \\
				ProxSVRG  & 0.1545      &  480    & 18573.3    \\
				RDA       & 0.1611      &  986    & 36001.3    \\
				SPStorm   & 0.1545      &  481    & 30068.4    \\
				SAGA      & \multicolumn{3}{c}{--}           \\
				\bottomrule
			\end{tabular}
		\end{subtable}
		\hfill
		\begin{subtable}{0.48\textwidth}
			\centering
			\caption{rcv1 dataset}
			\begin{tabular}{l S[table-format=1.6] S[table-format=4.0] S[table-format=5.1]}
				\toprule
				Algorithm & {$f$(best)} & {Iter.} & {Time (s)} \\
				\midrule
				{PSGA}      & {0.1264} & {64}    &   {219.2}     \\
				PStorm    & 0.1385      &  998    & 2066.6     \\
				ProxSVRG  & 0.1265      &  962    & 3055.2     \\
				RDA       & 0.1270      &  893    & 1207.8     \\
				SPStorm   & 0.1265      &  962    & 1997.6     \\
				SAGA      & 0.1265      &  961    & 72982.2    \\
				\bottomrule
			\end{tabular}
		\end{subtable}
	\end{table}
	
	In Table \ref{tab:convergence_all}, we observe that our algorithm(PSGA) obtains better objective function values $f(best)$. At the same time, our algorithm requires fewer iterations and less CPU time than other algorithms. Also we note that SAGA terminated immediately on the dataset news20  because the storage of the gradient look-up table exceeded the memory limit.
	
	\newpage
	\section{Conclusion}\label{Section5}
	In this paper, we propose a stochastic proximal gradient algorithm (PSGA) for solving composite convex optimization problems. Our method employs an adaptive step-size strategy, thereby relaxing both the strong convexity requirement of the objective function $f$ and fixed-step condition required by S-PStorm. In addition, our method employs an efficient variance reduction technique that reduces full gradient computation without requiring gradient storage. At the same time, we prove the gradient estimation error converges to zero almost surely. Moreover, we prove the strong convergence of our method and establish an $O(\sqrt{\frac{1}{k}})$ convergence rate. Numerical experiments on Logistic regression and Lasso regression illustrates the efficiency of our method.

	\section*{Declarations}
	\noindent{\bf Funding} The project is supported by National Natural Science Foundation of China No. 12101100 and Natural Science Foundation Project of Chongqing, China No. cstc2020jcyj-msxmX0738.\\
	\noindent{\bf Competing interests} The authors have no competing interests to declare that are relevant to the content of this article.\\

\end{document}